\newif\ifAMSvers
\AMSversfalse  
\newcommand\AMSonly[1]{\ifAMSvers{#1}\else\relax\fi}
\newcommand\notAMS[1]{\ifAMSvers{\relax}\else{\fi}}
\newcommand\AMSorNot[2]{\ifAMSvers{#1}\else{#2}\fi}
\ifAMSvers
  \documentclass[reqno]{ecgd-l}
\else
  \documentclass[twoside,12pt,letterpaper]{article}
  \usepackage[margin=1in,centering]{geometry}
\fi

\usepackage[singlelinecheck=false]{caption}
\usepackage{latexsym,graphicx}
\usepackage{graphbox}  
\usepackage{amssymb,amsmath,amsthm,amsopn,pifont} 
\ifAMSvers \usepackage{mathptmx}\fi
\usepackage[usenames, dvipsnames]{color}
\usepackage[implicit=true, 
colorlinks=true,linkcolor=Black,citecolor=Black,urlcolor=Black]{hyperref}%
\usepackage{aliascnt} 
\ifAMSvers\relax\else  \fi
\usepackage{longtable} 
\usepackage{array}     
\definecolor{light-gray}{gray}{0.8}
\newcolumntype{:}{!{\color{light-gray}\vline}} 

\usepackage{enumitem}
\setlist[enumerate,1]{label=\textbf{(\arabic*)},ref=(\arabic*)}
\usepackage{xifthen}
\ifAMSvers\relax\else
  \usepackage{placeins} 
  \newcommand\measurepage{\dimexpr\pagegoal-\pagetotal-\baselineskip\relax}
  \newtest{\IsThereSpaceOnPage}[1]{\lengthtest{\measurepage<#1}}
  \let\stdsection\section
  \renewcommand\section{\FloatBarrier%
    \ifthenelse{\IsThereSpaceOnPage{.15\textheight}}{\clearpage}{\relax}\stdsection}
  \let\stdsubsection\subsection
  \newcommand{\SubSecSkip}{\vspace{.5\baselineskip plus 8pt minus 3pt}}
  \renewcommand\subsection{\ifthenelse{\equal{\arabic{subsection}}{0}}%
   {\vspace{0pt plus 10pt}}{\SubSecSkip}\stdsubsection}
\fi 
\def\du#1{\underline{\underline{#1}}} 
\def\ssm{{\smallsetminus}}
\def\vecm{{\stackrel{\rightarrow}{m}}}
\def\vecx{{\stackrel{\rightarrow}{x}}}

\def\vecv{{\stackrel{\rightarrow}{v}}}

\newcommand{\mapstoself}%
  {\,\protect\rotatebox[origin=c]{90}{\large$\circlearrowleft$}\,}
\newcommand{\rmapsto}%
  {\mathrel{\raisebox{.22ex}{\protect\rotatebox[origin=c]{180}{$\mapsto$}}}}
\newcommand\ssk{{\smallskip}}
\newcommand\msk{{\medskip}}
\newcommand\bsk{{\bigskip}}

\newcommand\C{{\mathbb C}}
\newcommand\R{{\mathbb R}}

\def\Rhat{{\widehat\R}}

\def\x{{\bf x}}
\def\y{{\bf y}}
\newcommand\Z{{\mathbb Z}}
\usepackage{bm}
\newcommand{\dfn}[1]{{\everymath{\bm}\textbf{\textit{#1}}}}

\makeatletter
\def\newaliasedtheorem#1[#2]#3{%
  \newaliascnt{#1@alt}{#2}
  \newtheorem{#1}[#1@alt]{#3}
  \expandafter\newcommand\csname #1@altname\endcsname{#3}
}
\makeatother

\numberwithin{equation}{section}
\numberwithin{table}{section}

\newaliasedtheorem{lem}[equation]{Lemma}
\newaliasedtheorem{theo}[equation]{Theorem}
\newaliasedtheorem{cla}[equation]{Claim}
\newaliasedtheorem{coro}[equation]{Corollary}
\newaliasedtheorem{prop}[equation]{Proposition}
\newaliasedtheorem{question}[equation]{Question}
\newaliasedtheorem{conj}[equation]{Conjecture}
\newaliasedtheorem{asse}[equation]{Assertion}

\theoremstyle{definition}
\newaliasedtheorem{definition}[equation]{Definition}
\newaliasedtheorem{rem}[equation]{Remark}
\newaliasedtheorem{ex}[equation]{Example}

\usepackage[\AMSorNot{font={rm}}{font={sf}}]{caption}
\newcommand{\capf}{\relax}

\ifAMSvers\relax\else
  \newcommand\subjclass[2][2020]{
    \vspace{.2cm}\noindent%
    {\textbf{Mathematics Subject Classification (#1): }{#2}}\par}
  \newcommand\keywords[1]{\vspace{.2cm}\noindent\textbf{Keywords: }{#1\par}}
\fi

\begin{document}
\ifAMSvers
\else
  \thispagestyle{empty}
  \title{\textbf{The W. Thurston Algorithm\\Applied to Real Polynomial Maps }}
  \author{\textbf{Araceli Bonifant, John Milnor and Scott Sutherland}}
  \maketitle
  \date
\fi

\begin{abstract}
  This note will describe an effective procedure for constructing  critically
  finite real polynomial maps with specified combinatorics.
\end{abstract}


\section{Introduction.}\label{s-in}

The Thurston algorithm is a method for constructing critically
finite rational maps with specified combinatorics. (Compare \cite{DH}.)
In the general case, it requires quite a bit of work
even to describe the algorithm precisely, much less to prove convergence;
and the implementation is very difficult.

In this note we are concerned with the much easier case of real polynomial
maps with real critical points. In this case, the result can be stated quite
easily, and carried out without too much difficulty. (However the proof that
the algorithm always converges in the polynomial case depends essentially on
more difficult complex methods. We will simply refer to
Bielefeld-Fisher-Hubbard \cite{BFH} or Poirier \cite{P1} for this.)

\autoref{s1} will describe the data which must be presented
to the algorithm in order for it to produce a corresponding
 uniquely defined critically finite real polynomial map. The hardest step
 in carrying out the algorithm, at least when the degree is four or more, is
 the construction of polynomials with prescribed critical values.\footnote{
   The problem of understanding maps with specified critical values
   goes back to Hurwitz \cite{H}.}
\autoref{s2} will use methods  suggested by Douady and Sentenac to deal
with this problem. \autoref{s3} will then describe the actual
algorithm. \autoref{apCompute} deals with computational issues;
\autoref{s4} provides further examples,  in particular for
 the non-expansive case;  and \autoref{s5} is a brief
table providing more precise information about the various figures.
The authors plan to publish a sequel about real quadratic rational maps (see
\cite{BMS} for a preliminary version). 

Thurston's presentation was based on iteration in the Teichm\"uller space for
the Riemann sphere with finitely many marked points.
Many people have contributed to or applied
this theory or modifications of it. In addition to the papers cited
above, see especially \cite{Ba},  \cite{Ba-N}, \cite{HS}, and \cite{P2}, 
as well as \cite{HB-He}, \cite{Ch} , and \cite{J}.
{(The Hubbard-Schleicher  paper starts with an explicit description  of the
  Thurston algorithm for real quadratic polynomials.)} 
For a completely
different approach, see Dylan Thurston \cite{T}.\msk

For access to our maple code, and for an interactive demonstration program,
see\break \url{https://www.math.stonybrook.edu/~scott/ThurstonMethod/}
\bsk

\section{Combinatorics.}\label{s1}
Let $f:\R\to\R$ be a real polynomial map of degree $d\ge 2$
\AMSorNot{with}{which has} real critical
points {$c_1\le c_2\le\cdots\le c_{d-1}$}. The derivative of $f$ can
be written as
$$ f'(x) ~=~ a(x-c_1)\cdots(x-c_{d-1})~,\qquad{\rm with~leading~coefficient}
  \qquad a\ne 0~.$$
If we are given such an $f'\!$, then of course
$f$ is uniquely defined up to an additive constant.
By definition, the
{\dfn{real filled Julia set}} $K_\R(f)$ is the set of all real
numbers $x$ for which the forward orbit of $x$ is bounded.\msk

We will say that $f$ is in {\dfn{$K_\R\!$-normal form}}
if the smallest  point of its real filled Julia set is $x=0$, and the
largest one is $x=1$.
This is very convenient for graphical purposes, since it means  that all of the
interesting dynamics of $f$ can be observed by looking at its graph restricted
to the unit interval $I=[0,1]$, with all orbits outside of $I$ escaping to
infinity. Evidently $f$ can be put into $K_\R$-normal form by a unique
orientation preserving affine change of coordinates whenever
$K_\R(f)$ contains at least two distinct points.
(The exceptional cases where $K_\R(f)$ is empty, or consists of a single
non-attracting fixed point are of no interest to us.)

\begin{figure}[!htb]
\centerline{%
  \raisebox{-.5\height}{\includegraphics[width=.2\textwidth]{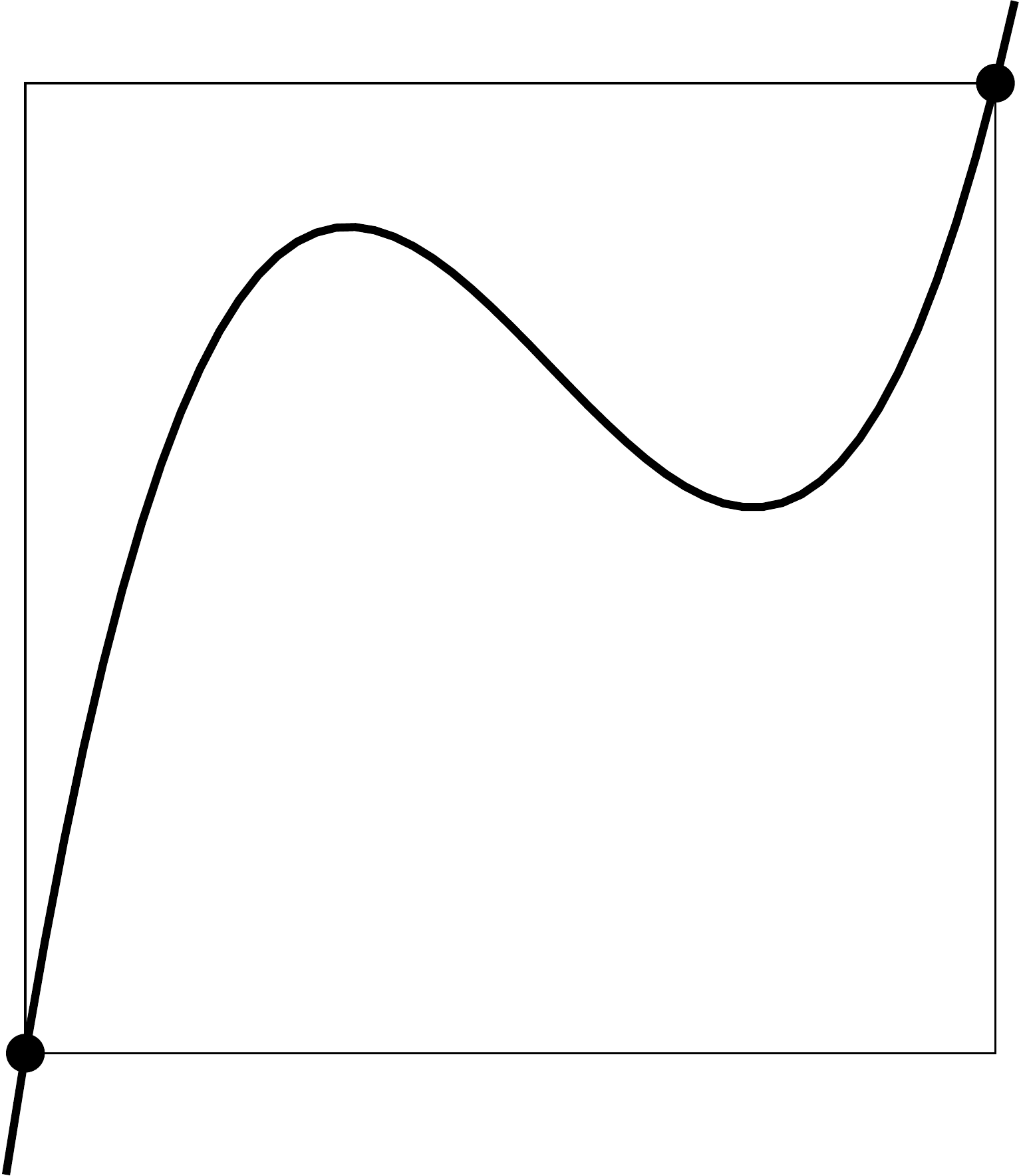}}\hfil
  \raisebox{-.5\height}{\includegraphics[width=.2\textwidth]{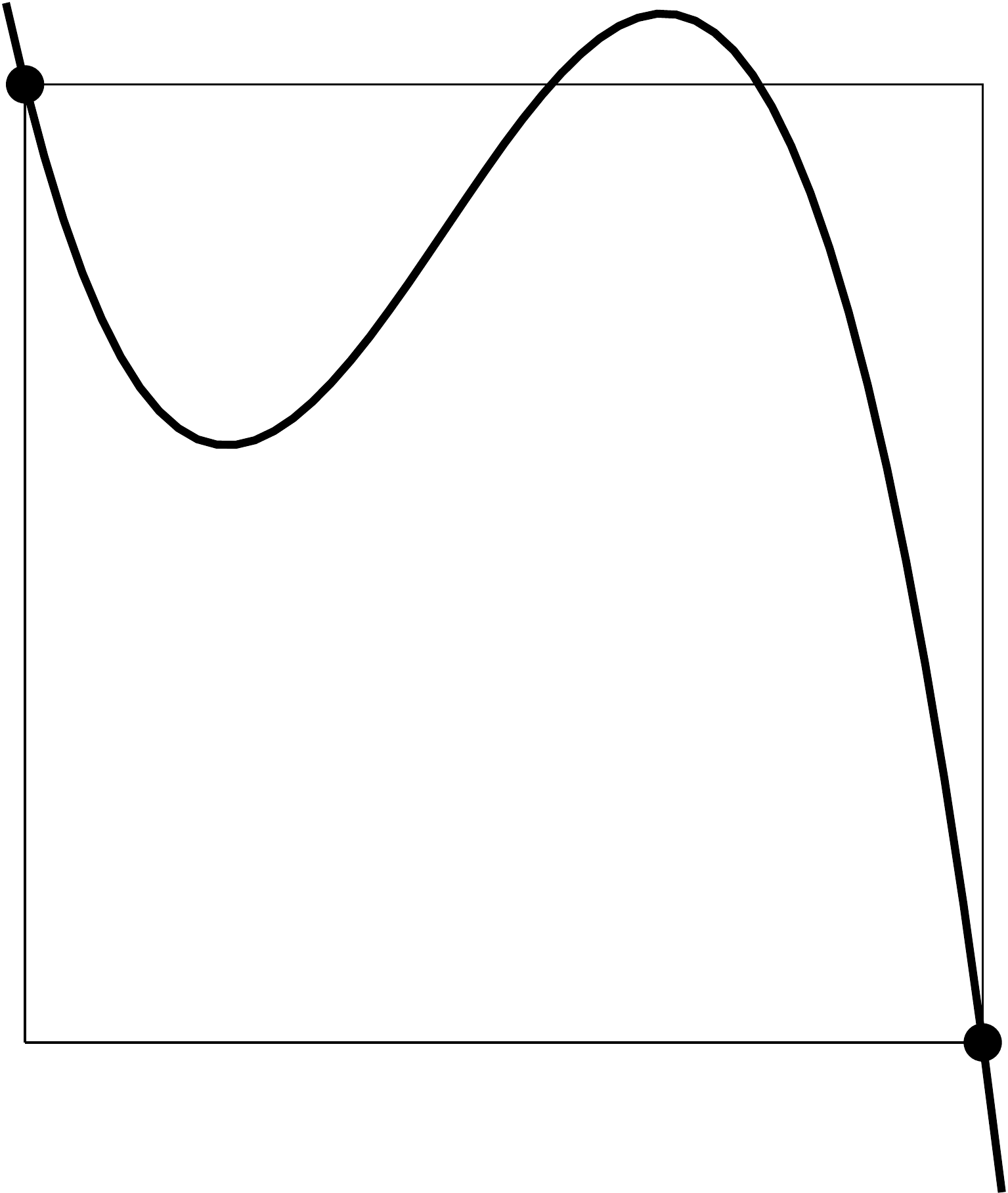}}\hfil
  \raisebox{-.45\height}{\includegraphics[width=.2\textwidth]{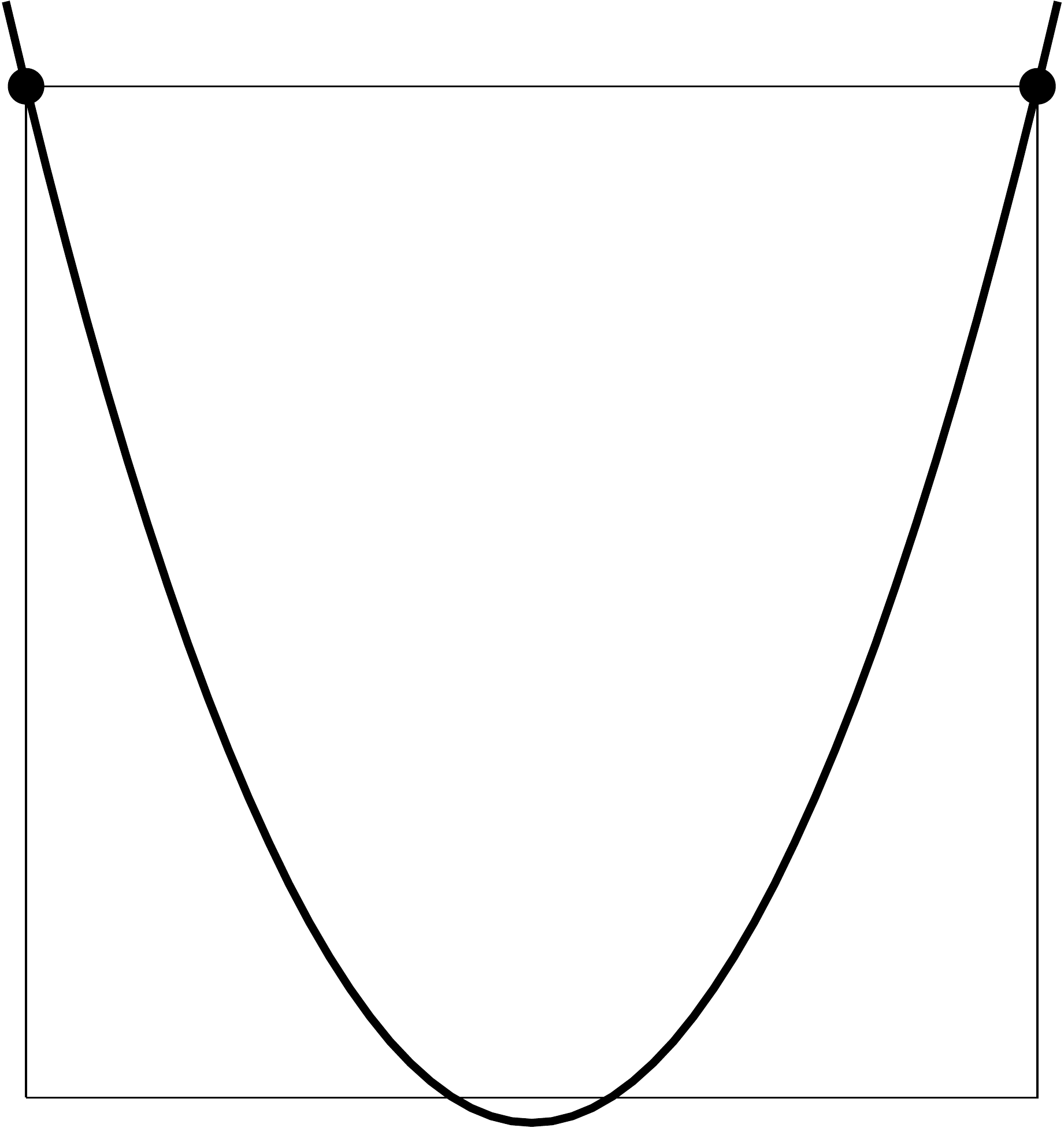}} \hfil
  \raisebox{-.5\height}{\includegraphics[width=.2\textwidth]{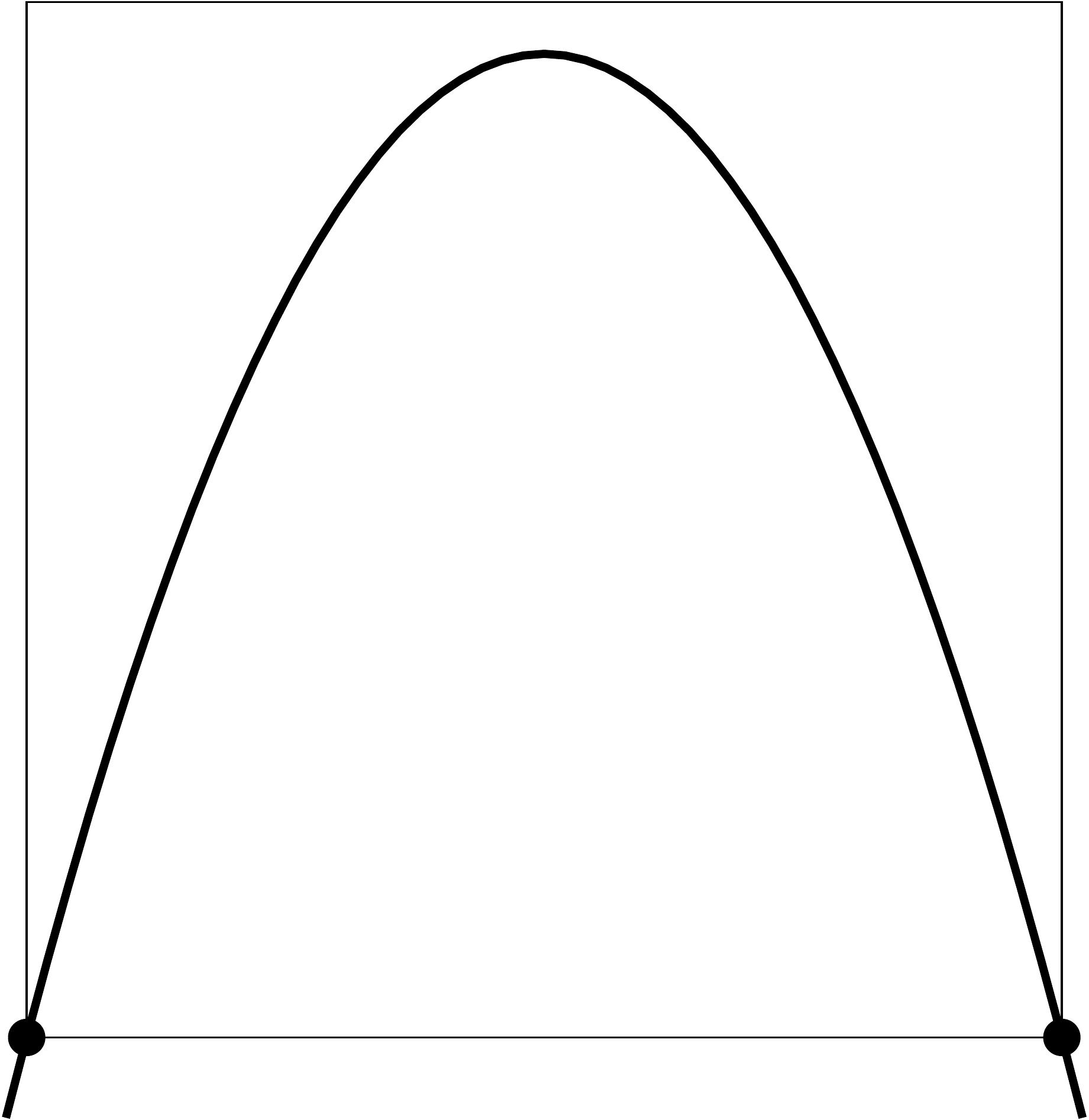}}
}

 \caption{\capf \label{f1} For a map of odd degree in $K_\R$-normal form, the end
 points $0$ and $1$ are either fixed points or form a period two orbit
 according as the leading coefficient is positive or negative. On the other
 hand, for maps of even degree, both end points map to one of the two. The
 middle two  graphs provide examples of maps for which $K_\R$ is not
 connected. In these examples, the critical points all lie in $[0,1]$
 but the critical values do not.}
\end{figure}

For any map in $K_\R$-normal form, note that the boundary $\{0,\,1\}$ of $I$
 necessarily maps into itself. There are four possible ways of
 mapping the boundary to itself, as illustrated in \autoref{f1}.
 \msk

\begin{definition}\label{def_combinatorics}
The map will be called {\dfn{critically finite}} if
every critical orbit is periodic or eventually periodic. A map of degree
$d$ has {\dfn{simple critical points}} if the $d-1$ critical
points are all distinct. That is, putting the map in $K_\R\!$ normal form,
they can be listed as
$$ 0<c_1<c_2<\cdots<c_{d-1}<1~.$$
(Later we will deal with higher order critical points, which may well be
boundary points of the interval $[0,1]$.  However, simple critical points are
always \dfn{interior} points of the interval.)
By the {\dfn{combinatorics}} of such a critically finite  map
with simple critical points
we will mean the sequence of integers $$\vecm~=~(m_0,~m_1,~\cdots~,~ m_n)$$
constructed as follows. Let
$$0=x_0<x_1<\cdots< x_n=1$$ be the list consisting of all points which are
either critical or postcritical, together with $0$~and~$1$ (if they are not
already included). Then there are unique integers $m_0, \cdots, m_n$
between zero and $n$ such that $f(x_j)= x_{m_j}$.\end{definition} 
\smallskip

Our goal is to show that 
the map $f$ is uniquely determined and effectively computable from its
combinatorics. More precisely, given a sequence
$~\vecm=(m_0,\cdots, m_n)$ satisfying a few simple necessary conditions,
there is one and only one corresponding polynomial map in $K_\R$
normal form, and the coefficients of this polynomial can be explicitly
computed.\smallskip

There will be a corresponding statement for the more general 
case where critical points of higher multiplicity are allowed;
and hence boundary critical points are also allowed.
However 
we will stick to the simple case for the moment.
\msk

\begin{figure}[!htb]
  \centerline{\includegraphics[width=.33\textwidth,align=c]{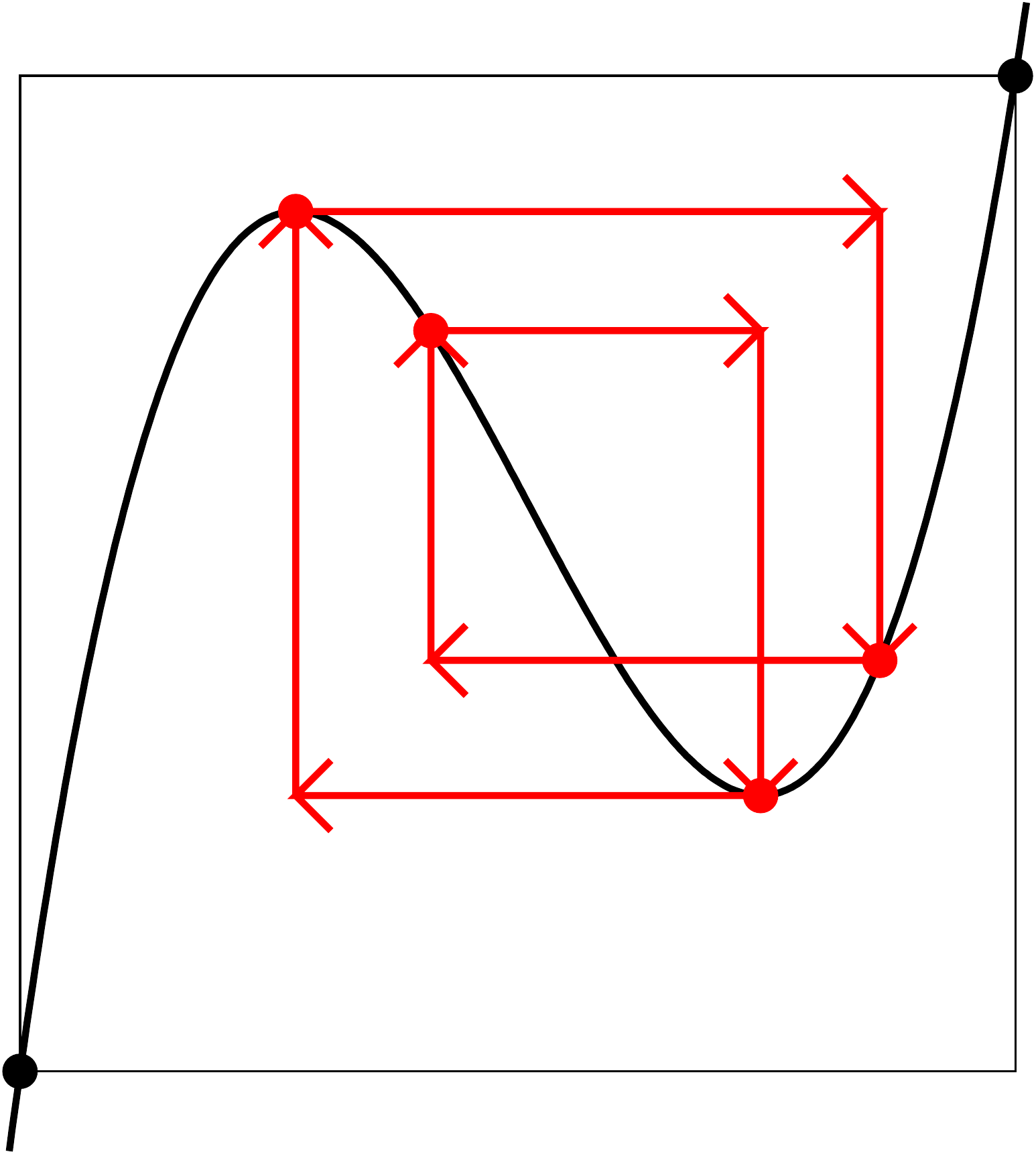}
 \hfil \includegraphics[width=.33\textwidth,align=c]{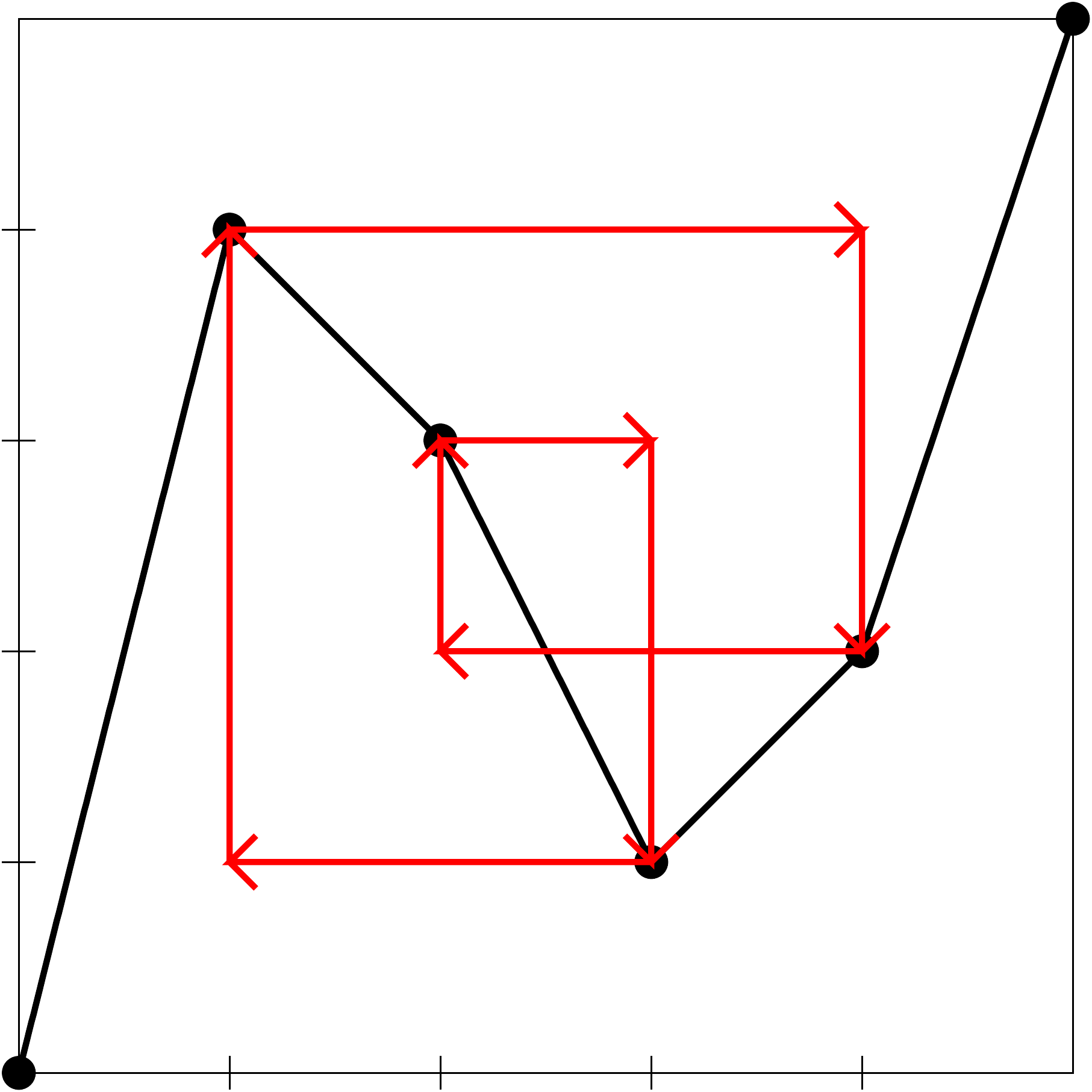}}

\caption{\capf \label{f2} On the left, example of a cubic map $f$ with
combinatorics $(0,4,3,1,2,5)$. Here the two critical points belong to a common
period four orbit. The red arrows describe the dynamics. On the right, a
corresponding graph for the piecewise-linear model of $f$, as described below.
See \autoref{coeff-tab} in \autoref{s5} for the precise equation 
of $f$. (For a neighborhood of $f$ in a complex parameter space
see \cite[Figure 13a]{BKM}.)}
\end{figure}

{\bf The piece-wise linear model.} 
It is often convenient to describe the 
combinatorics visually by considering
the graph of the 
function {$F\colon[0,\,n]\to [0,\,n]$}
 which maps each $j$ to $m_j$, and which is linear between
 integers.  Evidently the critical points of $f$ correspond
to the $d-1$ points $0<j<n$ where the graph has a local maximum or minimum.
(See \autoref{f2}.) In practice, we may replace $F$ by the rescaled map
$f(x)=F(n\,x)/n$  which sends the unit interval to itself.
\msk
 
Clearly the sequence $\vecm=(m_0,\cdots, m_n)$ must satisfy the following three 
restrictions:\ssk


\begin{enumerate}
\item $m_j\ne m_{j+1}$ for all $0\le j< n$. (Otherwise, the graph 
would have to be flat on the interval $~[m_j, m_{j+1}]$, or else have 
another $\max$ or $\min$ in the interior of this interval.)

\item The associated PL-graph must have at least one local maximum or
minimum with $0<j<n$, so that the degree satisfies $d\ge 2$.

\item {\bf(Framing)}\label{ComboFraming} $~~m_0$ must be equal to $0$ or~$n$,
and similarly $m_n$ must be equal to $0$ or~$n$.

\end{enumerate}\msk

\noindent In good cases, the following further condition will be satisfied.
See Bruin and Schleicher~\cite{BS} or Poirier \cite{P2}.  (The 
Bruin-Schleicher condition and the Poirier condition are stated differently;
but are  completely equivalent.) By a ``critical point'' of a piecewise linear
model map, we mean an interior local minimum or maximum. Note that the map
carries any edge (between successive integer points) onto a union of one
or more consecutive edges.

\begin{enumerate}[resume]
\item {\textbf{(Expansiveness\footnote{Caution: This terminology may be
confusing. The combinatorics is ``expansive'' if and only if the
associated pull-back transformation of \S\ref{s3} is contracting,
and hence convergent.})}}\label{ComboExpansive}
Every edge of the piecewise linear
model map must either have a critical boundary point, or else have an
iterated forward image which is long enough to contain
a critical point in its interior or boundary.
\end{enumerate}

\noindent 
Note that this last condition is always satisfied in the
  {\dfn{hyperbolic}} case, when every postcritical periodic cycle
  contains a critical point. For the behavior of examples which are not
  expansive, see \autoref{s4}.  In the limit map $f$, two or more
  of the $n+1$ points represented in the combinatorics will coalesce.
  
 One other condition is always satisfied for
   the combinatorics $\vecm$ constructed as above:\ssk

\begin{enumerate}[resume]
\item Each $j$ with $0<j<n$ is either a local maximum or minimum
for the associated PL-function, or is the iterated forward image of one.\ssk
\end{enumerate}

\noindent However this last condition is not really necessary. The proofs will
work just as well for choices of $\vecm$ which do not satisfy it.
\msk

\subsection*{Higher Order Critical Points.}
These definitions extend easily  to the case of polynomials with critical
points of 
higher multiplicity, and hence also allow
critical points on the boundary. In this case, each
of the integers $0\le j\le n$ must be assigned a
{\dfn{local degree}} 
$d_j\ge 1$ satisfying the following condition.

\begin{enumerate}[resume]
\item \label{localDegree} For all  $0<j<n$, the local degree $d_j$  
must be even if $j$ is a local minimum or maximum point for the associated
PL-graph, and odd otherwise. For the framing points $0$ and $n$, the
local degree must be one 
if the point is periodic, and odd in all cases. 
\end{enumerate}

\begin{figure}[ht!]
  \centerline{\includegraphics[width=.25\textwidth]{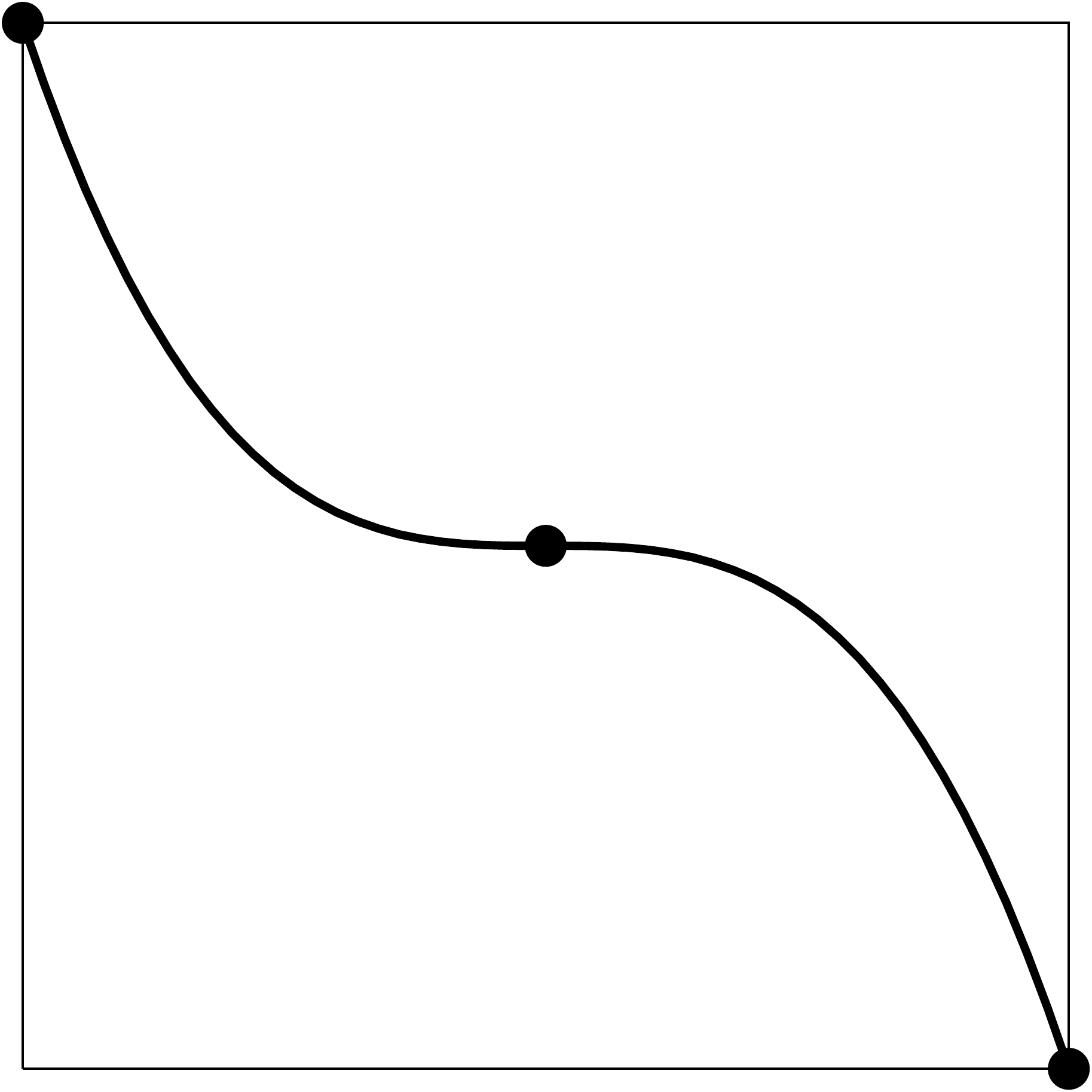}
      \hfil   \includegraphics[width=.25\textwidth]{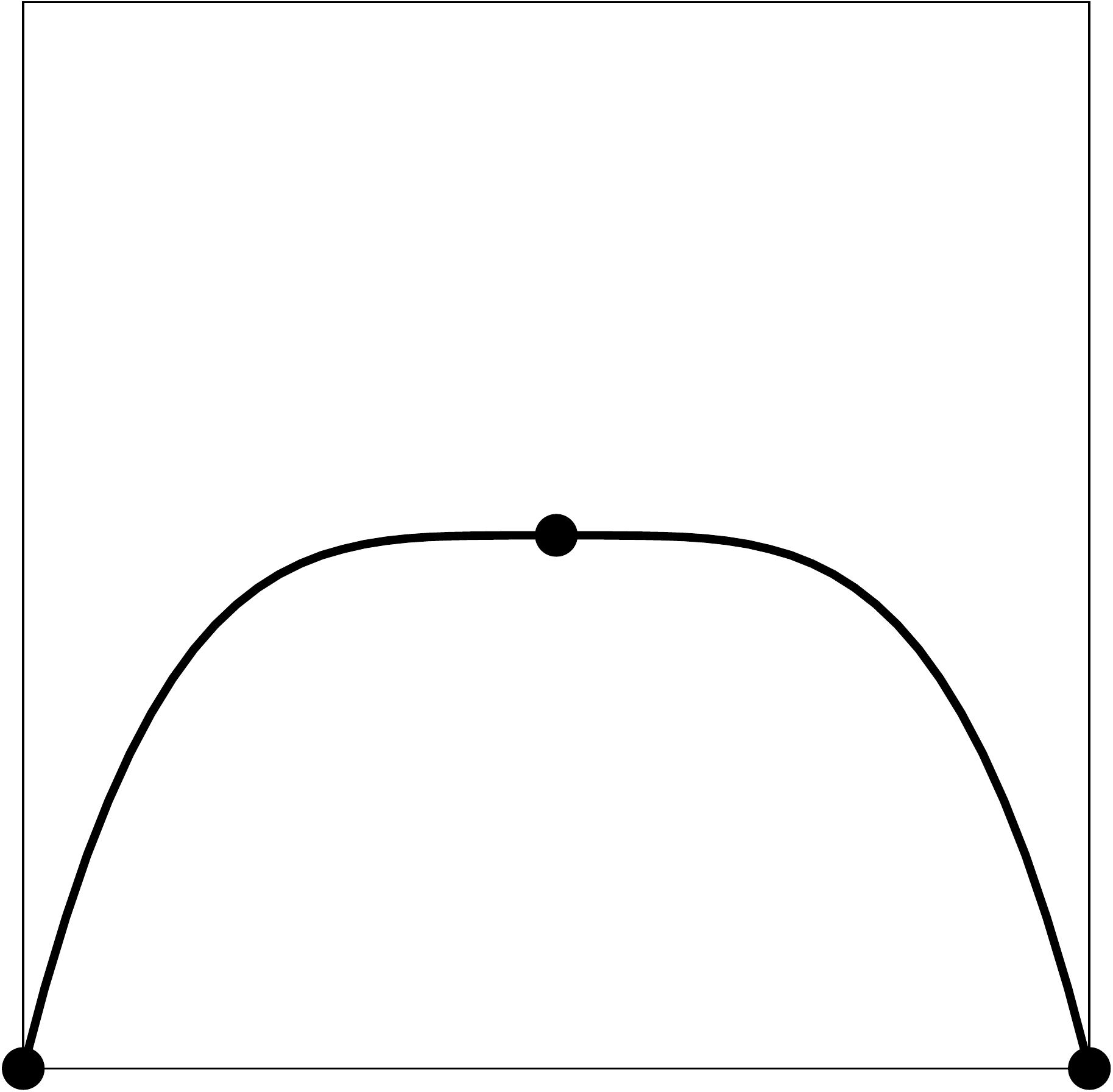}
      \hfil   \includegraphics[width=.25\textwidth]{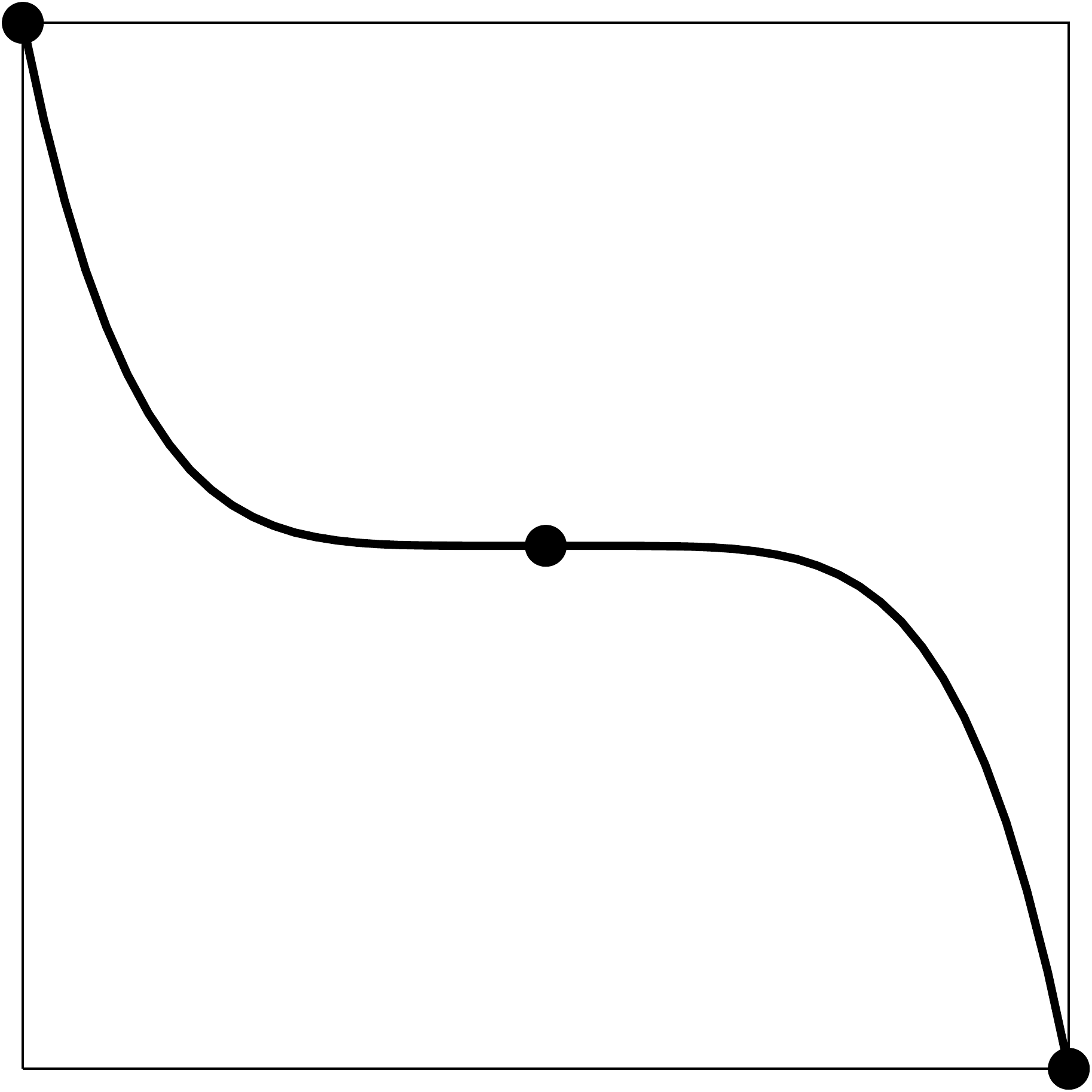}} 
    \caption{\capf \label{f3} Maps with equation
      $~ f(x)~=~  \big(1-(2x-1)^d\big)/2~$, for $d=3,4,5$. The
 combinatorics are respectively  $~~(2, 1^3,0)$,~~
      $(0,1^4,0)~$ and $~(2, 1^5, 0)$.}
\end{figure}

A convenient way of indicating this additional information
is to add $d_j$ as a formal superscript
on $m_j$ whenever $d_j>2$. (Compare \autoref{f3}.)\smallskip

\begin{figure}[!htb]
    \centerline{\includegraphics[width=.3\textwidth]{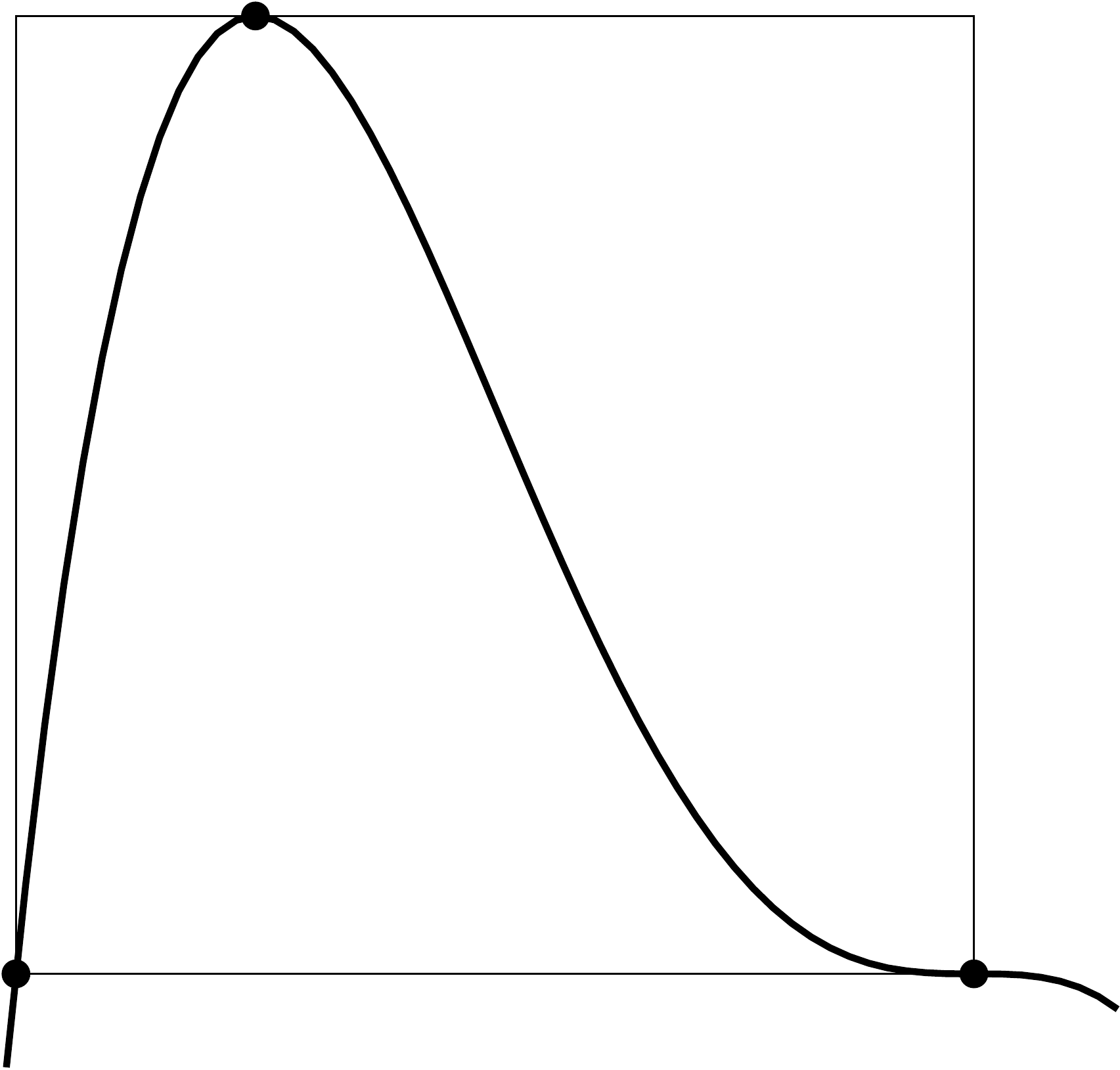}}
    \caption{\capf \label{f020} 
      The map $f(x)=kx(1-x)^3$ with $k=256/27$
      is critically finite with combinatorics $(0,2,0^3)$. Here the
      critical point is at $x=1/4$. The 
      mapping pattern   is
      $~~\du{x_1}\mapsto \underline{\du{x_2}} \mapsto x_0\mapstoself$.
      (
      The double or triple  underline indicates that the map has local
      degree two or three at the corresponding point.)}
      \end{figure}

To justify Condition~{\bf\ref{localDegree}} consider 
the corresponding polynomial map. If the degree for a fixed end 
point (or the product of degrees for a periodic cycle consisting of both
end points) were greater than one,
then the point or cycle would be attracting; hence nearby points
outside of $[0,1]$ would have bounded orbit, which is impossible.
Similarly, if a fixed end point had even local degree,
then all nearby points would map into $[0, 1]$, which again is impossible.
However, an endpoint which is not periodic can have an odd local
degree greater than one. (Compare \autoref{f020}.)
\ssk

Any point with $d_j>1$ will be referred to as a ``critical'' point. 
By definition, the associated total degree $d$ is the sum $1+\sum_j (d_j-1)$~.
\ssk

\autoref{f020}\ shows an example with
a critical point on the boundary. (In fact, all critical 
  values for \autoref{f020}\  
are also on the boundary of $f(\Rhat)$. Since there is only one framed
  polynomial map with the required critical value vector, it 
 follows that the Thurston
  algorithm for this combinatorics converges already on the first step.)
\ssk

\section{Prescribing Critical Values.}\label{s2}

This section will study the problem of finding a real polynomial map
with real critical points with a prescribed sequence of critical values.
The method of proof is due to Douady and
Sentenac. (Compare \cite[Appendix A]{MTr}.) \msk

First consider a polynomial $f$ of degree $d\ge 3$ with simple critical
points
$$ c_1<c_2<\cdots < c_{d-1}. $$
Let $(v_1,~v_2,~\cdots v_{d-1})$ be the corresponding sequence of
critical values $v_j=f(c_j)$.
\msk

\begin{theo}\label{sec3d-t1}
There exists a real polynomial with simple real critical points,
and with corresponding critical values
$(v_1,~v_2,~\cdots, ~v_{d-1})$ if and only if  the differences
$$ v_{j+1}-v_{j}\qquad{\rm for} \qquad 1\le j<d-1 $$
are all non-zero, and alternate in sign. The resulting polynomial 
$f$ is uniquely determined up to precomposition with an orientation
preserving affine change of variable, replacing $f(x)$ with $f(ax+b)$
where $a>0$. $($In particular, if there is a solution in $K_\R$
normal form, then it is uniquely determined.$)$
\end{theo}
\medskip

To begin the proof, the condition is clearly necessary, since  the polynomial
$f$ must be alternately monotone increasing or decreasing in the
intervals between critical points.
First consider the special case where we consider only maps such that
the derivative $f'$ is monic and centered, so that
$$ f'(x)~=~(x-c_1)\cdots(x-c_{d-1})\qquad
{\rm with}\qquad \sum_j c_j ~=~0~.$$
Let $v_j=f(c_j)$. 
Consider the sequences
$(\delta_1,\cdots, \delta_q)$ and $(s_1,~\cdots,~ s_q)$ of positive real
numbers, where $q=d-2$ and
$$ \delta_j = c_{j+1}-c_j~,\qquad s_j=|v_{j+1}-v_{j}|~.$$
Let $\R_+$ denote the set of all real numbers which are
$\ge 0$, and let $\R^q_+$ denote the $q$-fold product
$~\R_+\times\cdots\times\R_+$.\msk

\begin{lem}\label{sec3d-l1}
  There is a well defined map
$$ \Phi:(\delta_1,\cdots, \delta_q)\quad\mapsto\quad (s_1,\,\cdots,~s_q) $$
which sends  the interior of the space $\R_+^q$ diffeomorphically onto itself,
and extends to a map which sends each face of this product diffeomorphically
onto itself.
\end{lem}
\medskip

\begin{proof} Given the $\delta_j$, we can solve for each $c_j$ as a
linear function of $(\delta_1,\cdots,\delta_q)$ with constant coefficients.
Thus $f$ is uniquely determined up to an additive constant, hence the
differences $s_j$ are uniquely determined.

Next, given any $(\delta_1, \ldots, \delta_q)$ belonging to some face
of $\R^q_+$, we must show that the correspondence
$\Phi:(\delta_1, \ldots, \delta_q)\mapsto(s_1,\ldots, s_q)$ is locally a
diffeomorphism when restricted to that face.
Note that on a face of $\R^q_+$, the corresponding map will have non-simple critical
points. 
If $c_1\le\ldots\le c_{q+1}$ 
is the list of not-necessarily distinct critical points of the corresponding
polynomial $f$, then the distinct critical points can be listed\footnote{
  The case $r=1$ occurs when there is a single non-simple critical point,
  so all $s_j$ and $\delta_j$ are zero, and hence $\Phi$ is the identity map.}
as $\widehat{c}_1<\ldots<\widehat{c}_r$, with $2\le r\le q+1.$
The derivative
$g(x)=f'(x)$ can then be written as
$$ g(x)~=~(x-\widehat{c}_1)^{k_1}\cdots (x-\widehat{c}_r)^{k_r}~,$$
with multiplicities  $k_i\geq 1$ so that
$$k_1+\cdots+k_r=q+1\qquad {\rm and} \qquad k_1\widehat{c}_1+\cdots +
k_r\widehat{c}_r=~0~.$$
Each choice of the exponents  $~k_1,\cdots, k_r~$ corresponds
to the choice of some face of   the product  $\R_+^q$.\ssk
  
To deform $g$ within polynomials of this same form, set
$$\widehat{c}_i(t)  =\widehat{c}_i+ tw_i$$
where the $w_i$ are real numbers with $\sum k_iw_i =0$,
and where $t$ is a parameter which will tend to zero; and then set
$$g_t(x)~=~ \big(x-\widehat{c}_1(t)\big)^{k_1}\cdots
\big(x-\widehat{c}_r(t)\big)^{k_r}. $$
Since the logarithmic derivative of a product is equal to the sum of
logarithmic derivatives, and since $d\widehat{c}_i(t)/dt=w_i$,  we can write
$$ \frac{dg_t(x)/dt}{g_t(x)}~=~\sum_i \frac{-k_i w_i}{x-\widehat c_i(t)}~.$$
Multiplying both sides by $g_t(x)$ and then evaluating at $t=0$,
we see that the derivative $dg_t/dt$ at $t=0$ can be written as a product
$$~\big(x-\widehat{c}_1\big)^{k_1-1}\cdots\big(x-\widehat{c}_r\big)^{k_r-1}
h(x)~,$$ where
$$h(x)~=~ -\sum_i\left( k_i w_i  \prod_{\{j\,;\, j\neq i\}}
  (x-\widehat{c}_j)\right)~.$$
Here $h(x)$ is a non-zero
polynomial of degree at most $r-2$, since the coefficient
of the degree $r-1$  coefficient is $\sum -k_iw_i=0$.
Therefore the $r-1$  intervals $(\widehat c_i,~\widehat c_{i+1})$
cannot each contain a zero of $h(x)$; 
there must be at least one such interval
on which $h(x)$ has constant sign so that
the polynomial $dg_t/dt$  also has constant sign.

On the other hand, consider the function
\begin{equation*} t~\mapsto~
\widehat s_i(t)=\int_{\widehat{c}_i(t)}^{\widehat{c}_{i+1}(t)}|g_t(x)|\,dx~~=~~
\pm\int_{\widehat{c}_i(t)}^{\widehat{c}_{i+1}(t)}g_t(x)\,dx
\label{Es}~,   \end{equation*}
where $\pm$ is a sign which is determined by $i$ and $(k_1,\ldots,k_r)$.
Since $g_t(x)$ vanishes at the end points of the integration, we can
differentiate under the integral sign and then set $t=0$ to obtain
$$\left.\frac{d\widehat s_i(t)}{dt}\right|_{t=0}~=~
\pm\int_{\widehat c_i}^{\widehat c_{1+1}}\frac{dg_t}{dt}\, dx~.$$
This must be non-zero for at least one choice of $i$; which 
 proves that $\Phi$ restricted to any face of $\R^q_+$ is a local
diffeomorphism.\msk

It is not hard to check that $\Phi$ is a homogeneous map, with
$$\Phi(\lambda\delta_1,\cdots,\lambda\delta_q)
~=~\lambda^{q+2}\Phi(\delta_1,\cdots,\delta_q)~.$$

It follows easily that $\Phi$ carries each face to itself by a proper map,
which is necessarily a global diffeomorphism since each face is simply
connected. This proves the lemma. 
\end{proof}

As examples we have
\begin{description}
  \item[$q=1\AMSorNot{}{:}$] $\displaystyle \Phi(\delta_1) ~=~ \delta_1^3/2~$,
  \item[$q=2\AMSorNot{}{:}$] $\displaystyle
    \Phi(\delta_1,~\delta_2)~=~(\delta_1^4+ 2\delta_1^3\delta_2,~
    2\delta_1\delta_2^3+\delta_2^4)~.$

%

\end{description}

\begin{proof}[Proof of Theorem~\ref{sec3d-t1}]
Given $q+1$ critical values $v_j$, we can compute the
corresponding $s_j=|v_{j+1}-v_j|$ and hence find a polynomial $f(x)$
with positive leading coefficient which realizes this sequence
$(s_1,\cdots,s_q)$. There then exist unique coefficients $a\ne 0$ and $b$
so that the polynomial $~f(a x +b)~$ has the required critical values $v_j$.
\end{proof}
\smallskip

\begin{coro}\label{sec3d-c1}
Given $q+1$ prescribed critical values $v_j$ with
$v_{j+1}-v_j\ne 0$, there exists a polynomial in $K_\R$ normal form which
realizes these critical values if and only the differences alternate
in sign, and furthermore:

\begin{itemize}
\item either $v_1>0$ or $v_1<1$ according as $v_1>v_2$ or $v_1<v_2$, and

\item either $v_{q+1}>0$ or $v_{q+1}<1$ according as $v_{q+1}>v_q$ or
$v_{q+1}<v_q$.
\end{itemize}

\noindent The polynomial is unique when these conditions are satisfied.
\end{coro}

The proof is easily supplied. (Note that $K_\R$ will be equal to the unit
interval $[0,1]$ if and only if the $v_j$ all belong to $[0,1]$.)
\qed\bsk

\begin{rem}
Of course the Lemma~\ref{sec3d-l1} applies equally well to cases where there are
critical points with higher multiplicity. Correspondingly there are versions of
Theorem~\ref{sec3d-t1} 
and its Corollary~\ref{sec3d-c1} which apply to arbitrary combinatorics.
However the statements are somewhat more complicated.\end{rem} \bsk

\begin{figure}[htb!]
  \centerline{\includegraphics[width=.35\textwidth]{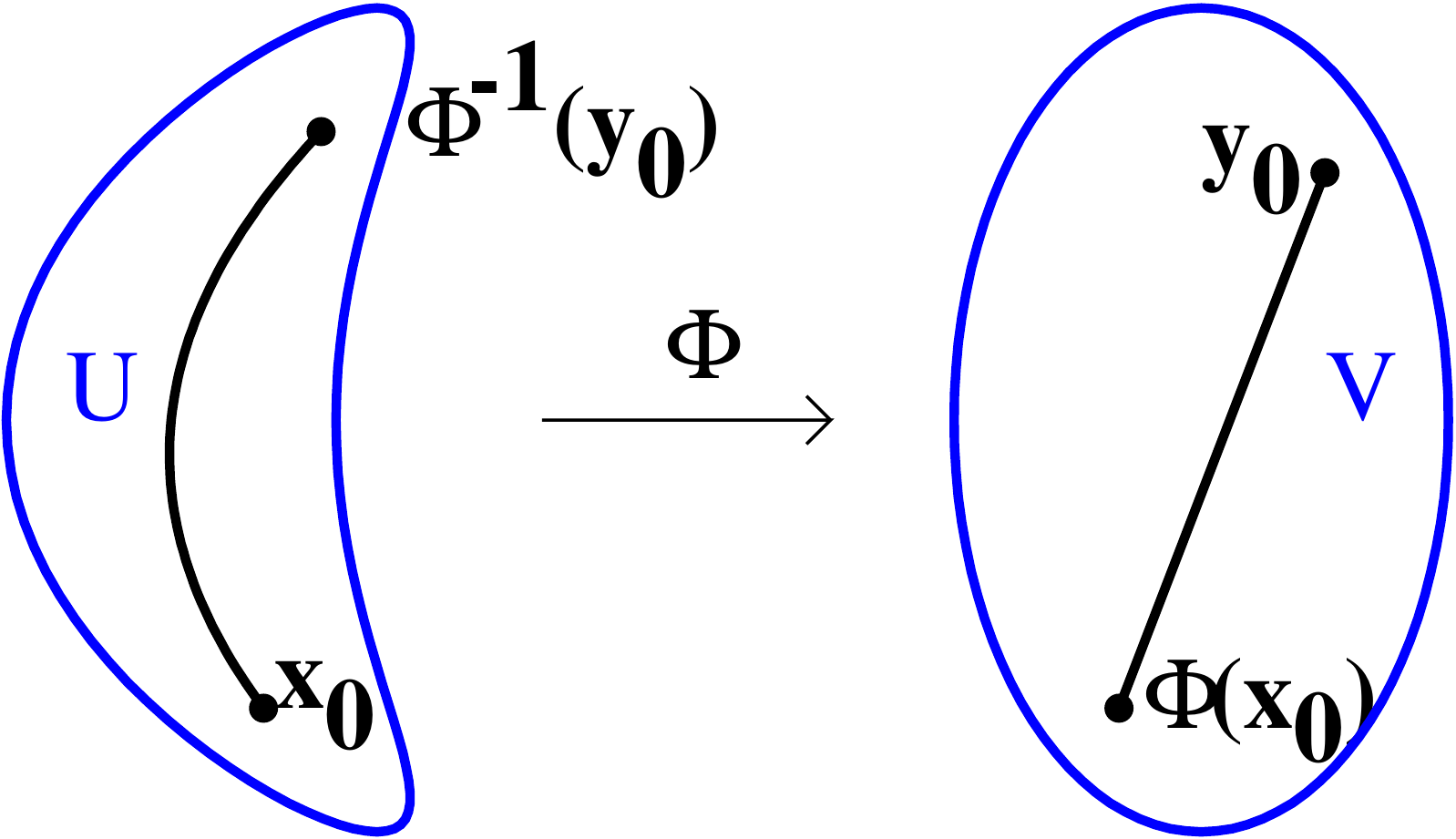}}
  \caption{\capf Computation of $\Phi^{-1}(\y_0)$.\label{sec3d-F1}}
\end{figure}

In order to actually carry out this construction, we must be able to
solve equations of the form
$$\Phi(\delta_1,\cdots,\delta_q)~=~(s_1,\cdots,s_q)$$
for any given sequence of $s_j$.

\begin{lem}\label{sec3d-L2} Let $\Phi:U\stackrel{\cong}{\longrightarrow} V$
  be an explicitly given diffeomorphism between open sets of $\R^n$ for any
  $n>0$, with $V$ convex. Then for any $\y_0\in V$ we can effectively
compute the pre-image $\Phi^{-1}(\y_0)\in U$.
\end{lem}\msk

\begin{proof} (Compare Figure~\ref{sec3d-F1}.)
Choose an arbitrary point $\x_0\in U$ and draw a straight line
$$ t ~\mapsto \y(t)~=~ (1-t)\Phi(\x_0)+t\,\y_0$$
from $\Phi(\x_0)$ to $\y_0$. Then the curve
$\x(t)=\Phi^{-1}\big(\y(t)\big)$ in $U$ joints $\x_0$ to the required
point $\Phi^{-1}(\y_0)$. This lifted curve satisfies a differential equation
of the form
$$ \Phi'(\x)\,\frac{d\x}{dt} ~~=~~ \frac{d\y}{dt} ~~=~~ \y_0-\Phi(\x_0)~.$$
Here $\Phi'$ is defined to be the $q\times q$ matrix of first partial
derivatives, and $d\x/dt$ and $d\y/dt$ should be thought of as column vectors.
Equivalently we can write
$$ \frac{d\x}{dt}~~=~~\Phi'(t)^{-1}\big(\y_0-\phi(\x_0)\big)~,$$
where the right hand side can be explicitly computed.
There exist standard packages for solving such systems of
differential equations to any reasonable degree of accuracy.
Therefore the curve $\x(t)$, and hence its endpoint $\Phi^{-1}(\y_0)$,
can be effectively computed.  
\end{proof}\msk

In fact, for the Douady-Sentenac diffeomorphism, a straightforward use of
Newton's method in several variables converges readily, \emph{provided} one
chooses an appropriate\footnote{See Appendix \ref{apCompute} for a more detailed
discussion of choice of starting point and a discussion of
implementation.} starting point.
For more general diffeomorphisms or poor choices of initial conditions,
Newton's method can behave very badly.
While the local convergence of Newton's method for such functions is well
understood and goes back over a century (see \cite{FineNewton}) and was greatly
expanded by Kantorovi\v{c} in the 1940s \cite{Kantorovich} (see also
\cite[Thm.~2.1]{Deuflhard}, \cite[\S2.8]{HubbardVect}),
the global behavior even for diffeomorphisms is not currently understood.



\msk
  
\section{The Algorithm.}\label{s3} Again we first consider the case of
simple critical points,  
and suppose that some
combinatorics $\vecm=(m_0,~m_1,~\cdots,~m_n)$ has been specified. Let
$X_n$ be the space \AMSorNot{consisting}{} of all $\vecx=(x_0,\cdots, x_n)$ with
$$0=x_0<x_1<\cdots<x_{n-1}<x_n=1~,$$
and let $0< j_1< \ldots <j_{d-1}<n~$   be the indices for
which $m_j$ is a local maximum or minimum.  The
{\dfn{pull-back transformation}} $$T:X_n\to X_n$$
which is associated with $\vecm$ can be described as follows. Given 
\AMSorNot{ any {$\vecx=(x_0,\cdots, x_n)$ in $X_n$},}
         { an arbitrary \hbox{$\vecx=(x_0,\cdots, x_n)$} in $X_n$,}
set $y_j=x_{m_j}$.  According to
\autoref{sec3d-c1}, there is a unique map 
  $f=f_\vecx$ in $K_\R$-normal form with critical values $v_i=y_{j_i}$.
  The image $T(\vecx)=\vecx'$ is a new element of $X_n$ satisfying the equation
  $$ f(x'_j)=y_j\qquad{\rm for~all}~~~j~.$$

\begin{figure}[!htb]
  \centerline{\includegraphics[height=.4\textwidth]{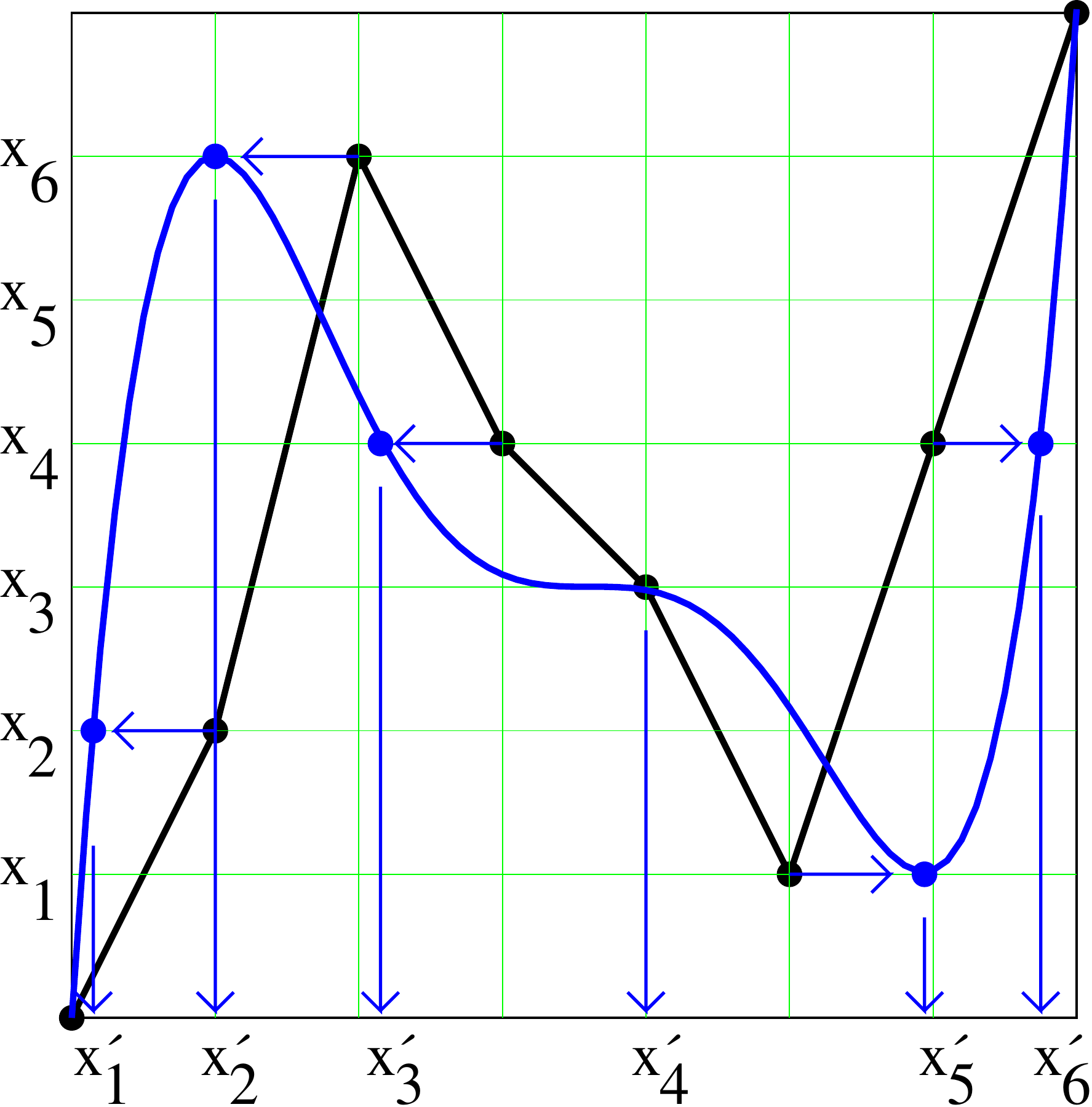}\hfil
              \includegraphics[height=.4\textwidth]{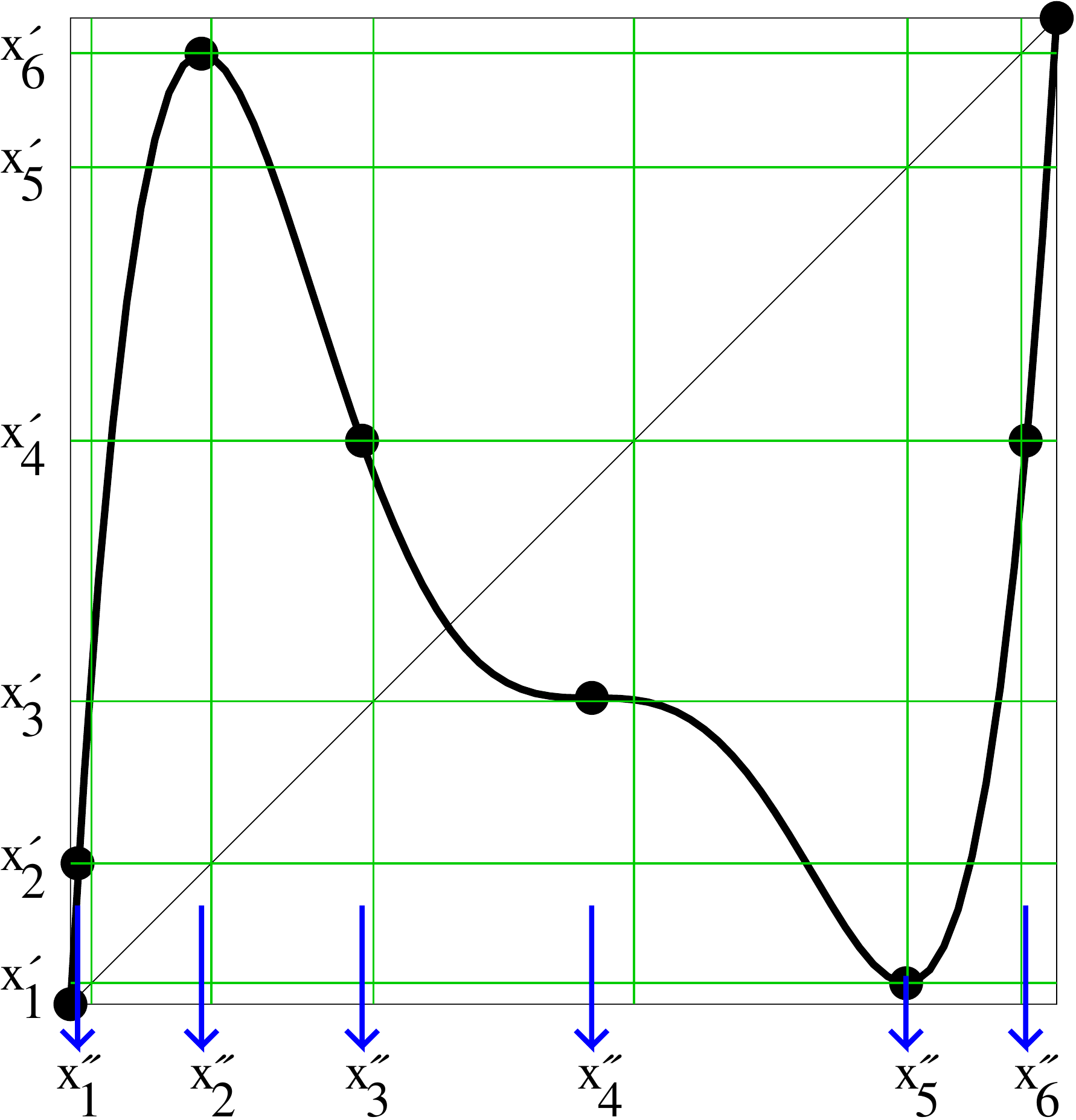}}
  \caption{\label{f-Thalg}\capf First two steps of the algorithm
    for the combinatorics $\vecm=(0,2,6^2,4,3^3,1^2,4,7)$.
    On the left:
    Start with the associated piecewise linear model map $f_0$, shown in black.
    The evenly spaced green grid line correspond to points $(x,y)$ with $x=x_j$
    or $y=x_j$. Note that the ``critical value vector'' for $f_0$ is the
    4-tuple $\vecv_0=(x_6,\,x_3,\,x_3,\,x_1)$, since the
    point  $(x_4',\,x_3)$ has been assigned local degree three.
 Let $f_1$, shown in blue,  be the unique framed polynomial map
    with the same critical value vector. Corresponding to each critical or
    postcritical point $(x, y)  =(x_j, x_{m_j})$
  on the graph of $f_0$, let $(x'_j,y_j)$ be the unique point on the
  corresponding lap of $f_1$ with the same $y$-coordinate. Note in particular
  that each critical point of $f_0$ corresponds to a critical point of $f_1$.
  The  correspondence
  $T:\vecx\mapsto\vecx'$ is called the  pull-back transformation. On the right:
  The green grid coordinates are now the points $(x,y)$ with $x=x'_j$ or $y=x'_j$.
  Let $f_2$ be the polynomial map with critical value vector
  $(x'_6, x'_3, x'_3, x'_1)$. Now for each $(x'_j, x'_{m_j})$ let
  $(x''_j,\,x'_{m_j})$ be the point on the corresponding lap of $f_2$
  with height $x'_{m_j}$. Then the pull-back transformation sends
  $\vecx'$ to $\vecx''$. Note that each $x''_j$ is fairly close to $x'_j$,
so that the iteration seems to be converging well.}
\end{figure}

\begin{figure}[!htb]
  \centerline{\includegraphics[width=.4\textwidth]{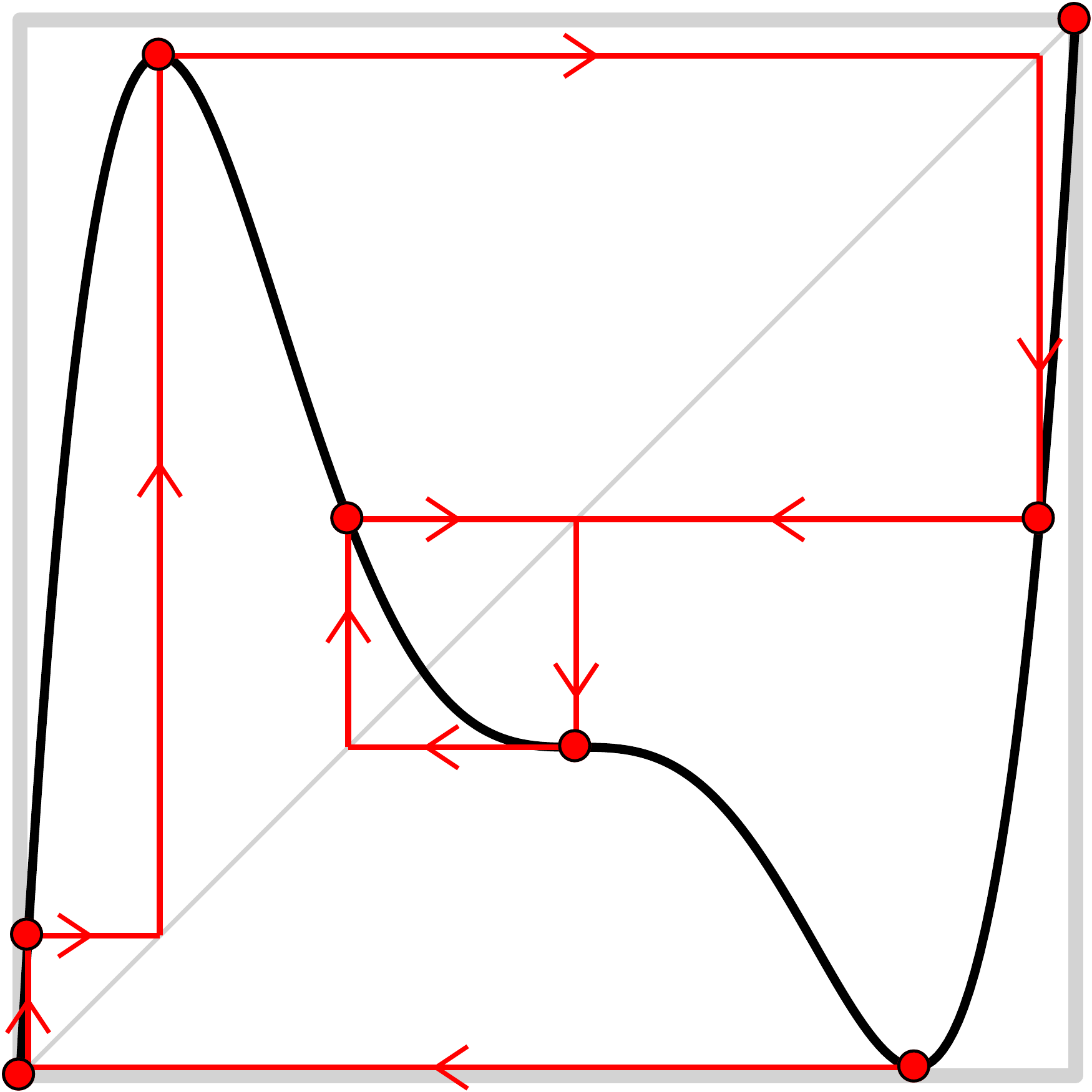}} 
  \caption{\label{f-thex-lim} \capf Showing the limiting map for the combinatorics
    $(0,2,6^2,4,3^3,1^2,4,7)$ of \autoref{f-Thalg}. (Here the last critical
    value is strictly positive, although very close to zero.) The mapping pattern
is  $~\du{x_5}\mapsto x_1\mapsto\du{x_2}\mapsto  x_6\mapsto x_4\leftrightarrow\underline{\du{x_3}}$.}
\end{figure}

\ssk
To see that $\vecx'$ is uniquely defined, first consider the $d-1$ critical 
indices~$j_i$. Choose $x'_{j_i}$ to be the corresponding critical point of $f$,
so that   $$f(x'_{j_i})={y_{j_i}}=v_i~,\qquad{\rm as~required.}$$
Now consider the remaining indices $j$. The $d-1$ critical points of $f$
separate the graph of $f$ into $d$ monotone segments called laps.\footnote{
   By definition, a {\dfn{lap}} of a piecewise monotone function
   is a maximal interval of monotonicity. In the special case of maps
   with simple critical points, this is the same as an interval between
   consecutive critical points or end points.}
Since we
require that  $x'_0<x'_1<\cdots<x'_n$, each of the remaining $x'_j$ must
belong to a well defined lap.  Therefore, using the intermediate value theorem, 
$x'_j$ is uniquely determined by its image $f(x'_j)=y_j$
 within that lap. This completes the description of $T$ in the special case
 of simple critical points.
 A very similar argument applies in the more general case, as long as
 Condition~\ref{localDegree} of \autoref{s2} holds; see \autoref{f-Thalg}.
 See \autoref{apCompute} for a more detailed discussion.

\msk

 If we can find a fixed point $\vecx=T(\vecx)$ for this transformation, then
 evidently the associated $f$ will be the desired critically finite 
 polynomial which satisfies the identity
 \begin{equation*}\label{E-fp} f(x_j)~=~x_{m_j}~. \end{equation*}
 In particular, starting with any $\vecx\in X_n$, if the successive
 iterated images $T^{\circ k}(\vecx)$ converge to a limit in $X_n$, then
 the polynomial  $f_\infty$ associated with this limit will be
 such a fixed point. 
 In fact if the combinatorics satisfies all of our
 requirements (including the expansiveness condition), then      
 such a unique  limit always exists. 
  A more detailed  explanation follows.

\bsk


\begin{definition}\label{d-Thurstonmap} 
  A {\dfn{Thurston map}} is an orientation preserving branched
  covering map from a topological 2-sphere onto itself which is
  ``critically finite'' in the sense that every branch point has a finite
  forward orbit.
  (It will be convenient to refer to the branch points as  critical points.)
  If there are $n\ge 5$ postcritical points (and in many cases if there are $n=4$),
 Thurston defines the pull-back transformation
 on an associated Teichm\"uller space for the surface of genus zero with $n$
 marked points, and proves that this transformation converges to a unique limit
 unless there is a well defined obstruction which prevents convergence.
 See \cite{DH} for details.\ssk

 As examples, flexible Latt\`es maps (see e.g. \cite{M}) 
provide 1-parameter families of
rational  maps of arbitrarily high degree with 4 postcritical points,
 all in the same equivalence class
 as Thurston maps. Such examples cannot occur when there are 5 or more
 postcritical points. \ssk

  Such a Thurston map is a {\dfn{topological
polynomial}} if there is a marked branch point (corresponding to the point at
infinity for an actual polynomial) which is fixed, and has no preimages other
than itself. This special case is much easier to deal with. (Compare
\cite{BFH}.)
In particular, the only possible obstruction is a Levy
cycle\footnote{ A Levy cycle is a special kind of
Thurston obstruction. A {\dfn{multicurve}}
$\Gamma=\{\gamma_1, \ldots, \gamma_n\}$ for the map $f$
is defined as a collection of disjoint, non-homotopic, and non-peripheral
simple closed curves in $\widehat\C\ssm P_f$, where $P_f$ is the
postcritical set. (Compare \cite{DH}.)
These form a {\dfn{Levy cycle}} if the curves
can be numbered so that, for each $i$ modulo $n$ the curve
$\gamma_{i-1}$ is homotopic to some component $\gamma_i'$ of $f^{-1}(\gamma_i)$
in $\widehat\C\ssm P_f$, such that
the map $f:\gamma_i' \rightarrow \gamma_i$ is a 
homeomorphism. (Compare \cite[Definition 5.2]{BFH}.)};
the cases with $n\le 4$  postcritical points present no problem.
Such a topological polynomial can always be represented by a Hubbard tree,
and there is no Levy
cycle if and only if an appropriate expansiveness condition is satisfied.
(Compare \cite{BS} or \cite{P2}.)\ssk

By a {\dfn{real topological polynomial}}
we will mean a piecewise monotone
map from $\R/\Z$ to itself such that only zero maps to zero, with specified
``critical points'' satisfying the requirements of Condition \ref{localDegree}
of \autoref{s2}. This is much easier to deal with than a complex
topological polynomial.
Each choice of combinatorics (whether or not the Expansiveness Condition
\ref{ComboExpansive}
of \autoref{s1} is satisfied) gives rise to a formal Hubbard tree contained
in the real line, with its associated real topological
polynomial. The Thurston algorithm for this real topological
polynomial will converge to an actual polynomial having the specified 
combinatorics if and only if the expansiveness condition
is satisfied.
(Actually, we will see in \autoref{s4} that there is a weaker form of
convergence even without expansiveness.)
\end{definition}

Of course we need an actual complex topological polynomial in order
to apply Thurston's convergence theorem. 
We will describe the construction of the complex topological polynomial from the
real combinatorics in one typical example, leaving the general construction to
the reader.

      
\begin{figure}[!htb]
  \centerline{\includegraphics[width=\AMSorNot{\textwidth}{.9\textwidth}]{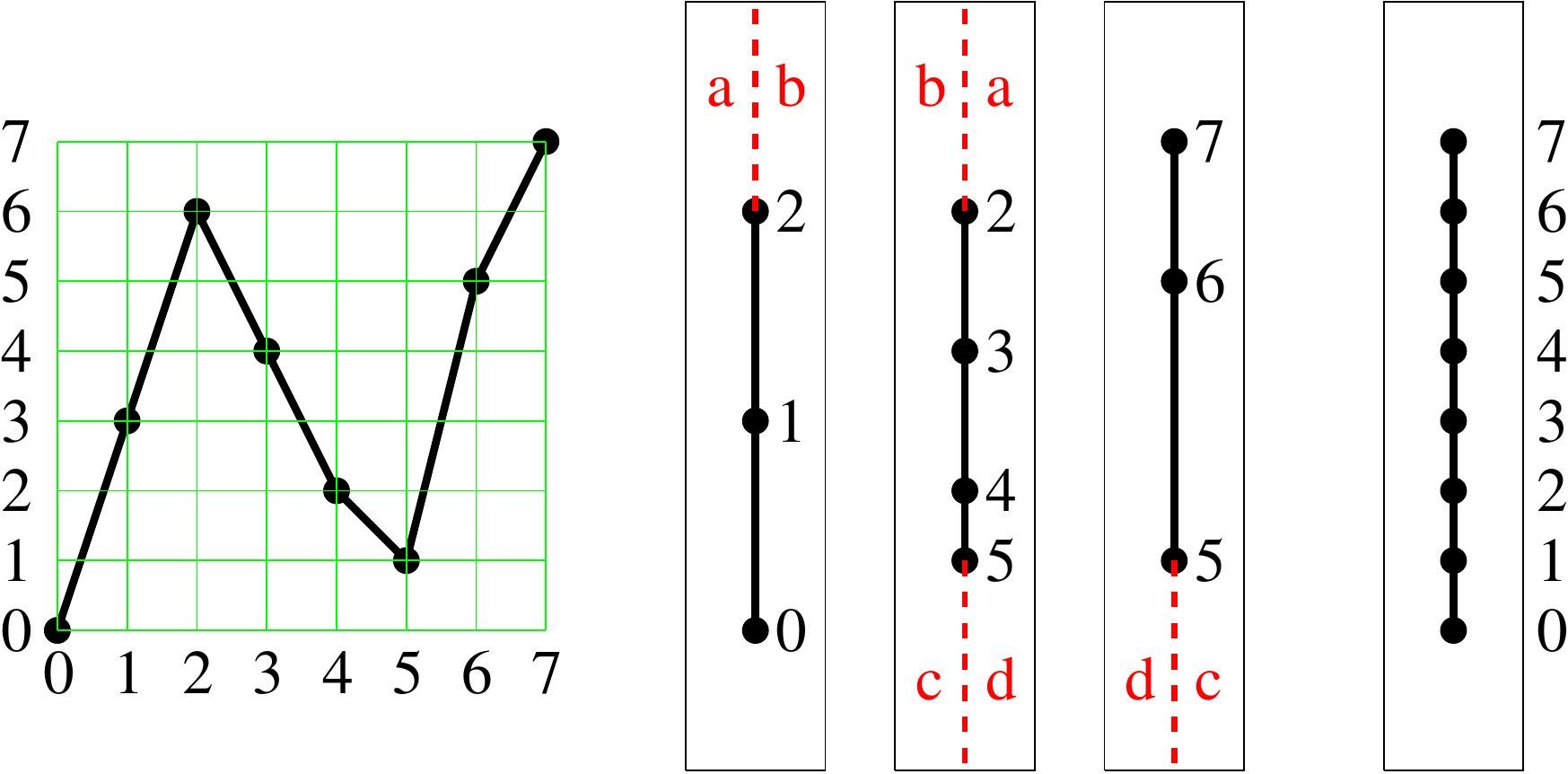}}
  \caption{\label{f-toppoly} \capf
    Construction of a branched covering
    of the sphere with combinatorics $(0,3,6,4,2,1,5,7)$.}
  \end{figure}
  
  Start with the piecewise linear model (\autoref{f-toppoly}-left),
  with   combinatorics 
$$(0,\,3,\,6,\,4,\,2,\,1,\,5,\,7)~.$$

{\bf Step 1.}  For each of the three laps, choose a copy of
$\C$ with the real axis vertical, as indicated schematically in the middle
of the figure, and project each lap  to the corresponding real axis.
For each marked point on the lap, mark a corresponding point
on this real axis, with height the associated $y$ value, but with label
the associated $x$ value.

Next slit each of these copies of $\C$ along the dotted lines, and then
paste the resulting boundary curves together in pairs;
pasting side ``a'' to side ``a'' and so on. 
The result will be a connected
simply connected surface $S$, which is homeomorphic to $\C$.\ssk

{\bf Step 2.} Project $S$
onto another copy $S'$ of $\C$, as represented on the right of the
figure, where now the marked points are labeled by their $y$ values.
Note that this projection is a branched covering, locally two-to-one
at the branch points, which are labeled {\bf 2} and {\bf 5}. 
Thus we have a branched covering  $S\to S'$.\ssk

{\bf Step 3.} Finally identify $S$ with $S'$ by choosing
a homeomorphism $S\leftrightarrow S'$  which sends the emphasized
part of each real axis to the real axis, and sends
each marked point to the point with the same label. Thus 
(after adding a point at infinity to each surface)
we obtain the required map from a topological sphere to itself.
As an example, the point labeled {\bf 1} on $S'$ is identified
with the point labeled {\bf 1} on $S$, and maps to the point
labeled {\bf 3} on $S'$.\ssk

The discussion would be similar for any example with distinct real critical
points.  The case with higher order critical points is more complicated, and
can be dealt with more easily by using external rays in the complex plane.
 See for example \cite{P1} and \cite{P2}.\bsk

\appendix
\section{Implementing the method computationally}\label{apCompute}

Implementation of the method is relatively straightforward, although there are
a few issues which need a little care. For low degree polynomials with fairly
tame combinatorics, all of the calculations can be done in standard
double-precision. For polynomials of degree~6 or higher (and in some
particular low-degree cases), calculations often require more digits in order
to converge.

\begin{figure}
  \centerline{\includegraphics[width=.35\textwidth]{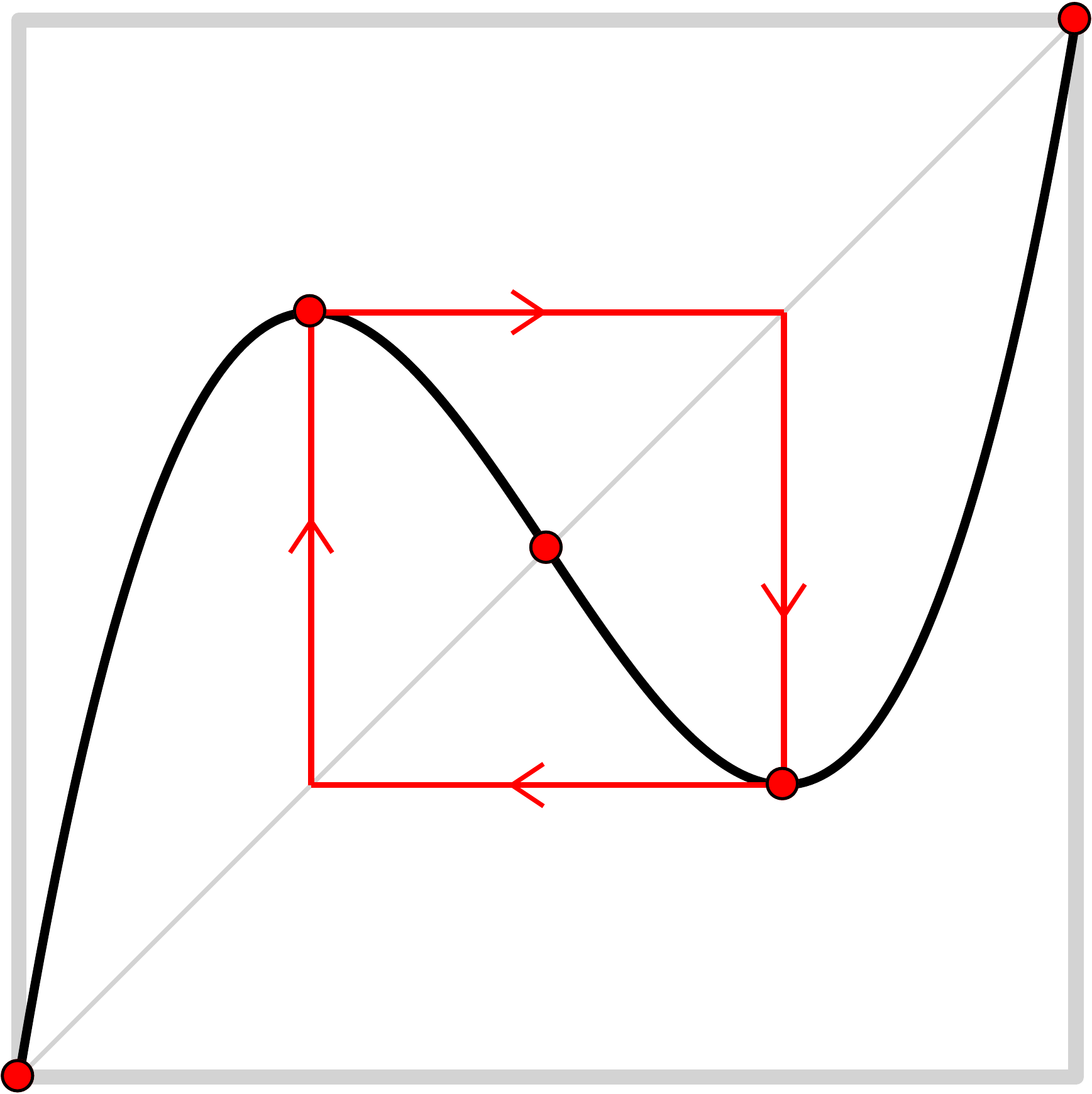} \hfil
    \includegraphics[width=.35\textwidth]{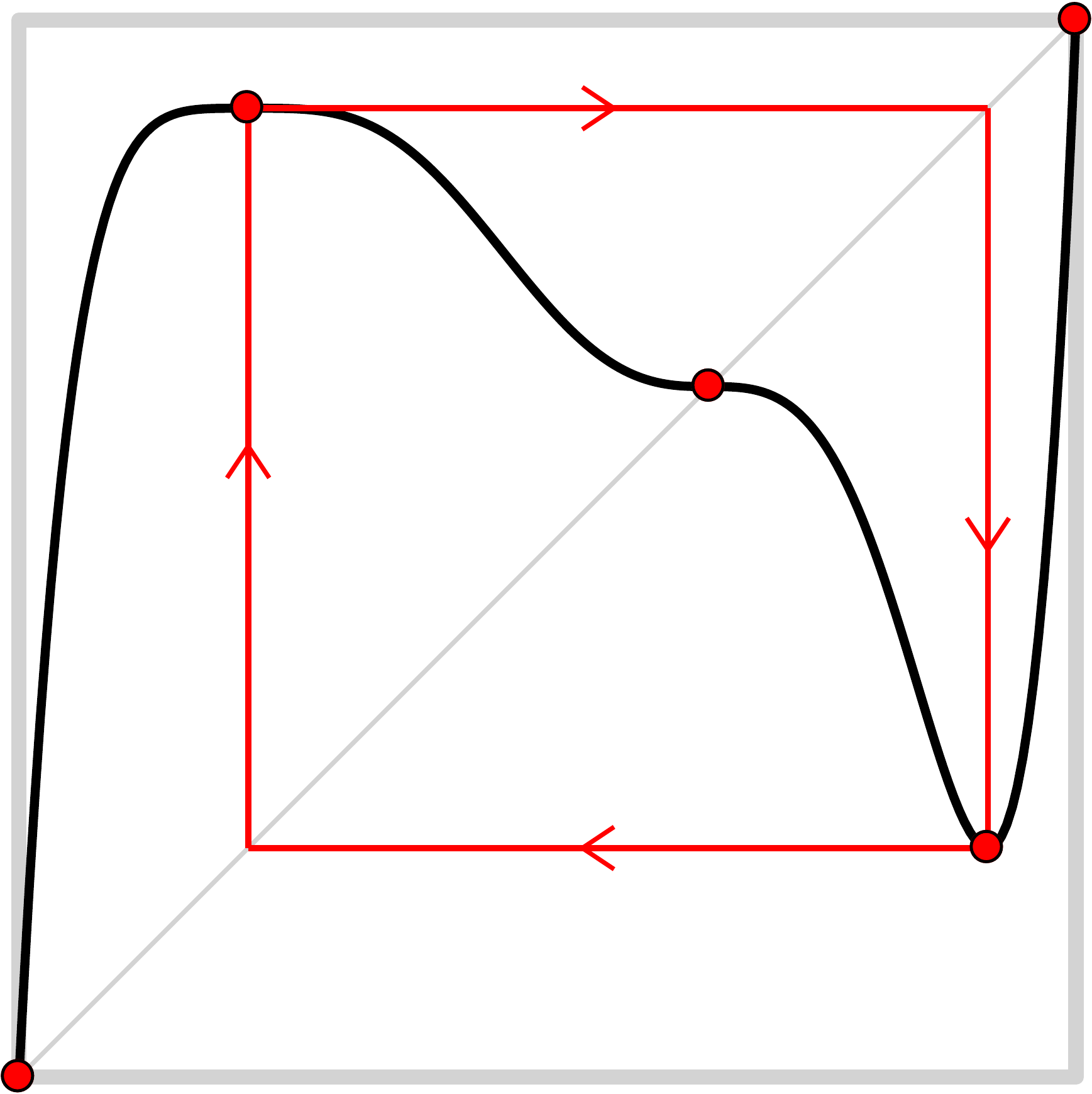}  }
\caption{\label{f-03214} \capf
 Two limiting polynomials, both with topological combinatorics
   $(0,3,2,1,4)$ (that is, the same piecewise-linear map).
  The polynomial on the left is cubic with simple critical points, while the
  polynomial on the right has degree~7 and combinatorics $(0,3^4,2^3,1,4)$.
}
\end{figure}

As an explicit example, for the combinatorics  $(0,3,2,1,4)$ the method
converges quite rapidly to the degree~3 limiting polynomial; see
\autoref{f-03214}. Using double precision (about~13 decimal digits), the method
converges to within $6\times10^{-6}$ in 5~steps, or better than $10^{-12}$
after 17~steps. 
For comparison, if we keep the combinatorics the same but change the first
critical point to have local degree~4 and make the central fixed point have
local degree~3 (that is, combinatorics $(0,3^4,2^3,1,4)$), the method requires
much more precision to converge to the corresponding degree~7 polynomial. Using 
double-precision arithmetic, it converges to within $1.3\times10^{-5}$ after
8~steps, but then loses precision, with the error oscillating between
$2\times10^{-4}$ and $7\times10^{-6}$.  Increasing the precision to
20~digits gives better than $10^{-7}$ in 10~steps.  One needs at least
27~digits of precision to get a limiting map good to within $10^{-12}$ (in
19~steps).  

We have written an implementation\footnote{%
see \url{https://www.math.stonybrook.edu/~scott/ThurstonMethod/}}
for arbitrary combinatorics and degree in
\textsf{Maple}, although it would be straightforward to port this to many
other languages, as long as the language supports multiple-precision
arithmetic. 

\bsk
\subsection*{Implementation}
Naturally, it is important to begin with combinatorics that are topologically
possible and fully describe the situation, as described in \autoref{s1}.

We insist that the framing points ($x=0$ and $x=1$ for a map in
$K_\R$ form) must be specified as part of the combinatorics.  For polynomials
of odd degree, we assume that the two end points are either fixed or form a
period~2 cycle, while for even degree we assume $x=0$ is fixed and $x=1$ is
the preimage of $x=0$.  The other even degree case is easily obtained from
this by a change of coordinates.
\bsk

\noindent
In the discussion below, we will use $x_m^{[\ell]}$ 
to denote the location of the $m$th marked point at the $\ell$th step of
the process. We will also use $c_j^{[\ell]}$ to denote the $j$th critical point
(which is, of course, one of the $x_m^{[\ell]}$).
When the particular step is irrelevant or apparent, we may omit the
superscript.  

\medbreak
We perform the following steps:
\begin{enumerate}[start=0, ref={\textbf{(\arabic*)}}]
\item \label{initstep} \textbf{Initialization.}
From the combinatorics (see \autoref{def_combinatorics}), the indices
corresponding to critical  
points can be inferred if local degrees are not explicitly specified.  
(As noted in \S\ref{s1}, 
when local degrees are not specified, we assume that all
critical points have local degree~2 and all other points are regular
points.)  Since the critical points are not required to be simple, we will
have $r$ distinct\footnote{Since the critical points here must be
  distinct but not necessarily simple, we have simplified the notation from
  \autoref{sec3d-l1} and write $c_j$ instead of $\widehat{c}_j$.}
  critical points $(c_1, c_2, \ldots, c_r)$ and $r$ (not
necessarily distinct) critical values $(v_1,\ldots,v_r)$.

%
Further, the laps can be determined just by inspection of the
given combinatorics, with laps bounded by each turning point.  Critical
points of odd degree can be ignored when defining laps, because they
do not affect the covering properties of the map.
Here, because we will need to solve numerically for
the two framing points, it is important to temporarily
allow the initial and final 
laps to extend sufficiently far\footnote{For example, for the
  combinatorics $(0,2,0^3)$ of \autoref{f020}, the right-hand
  framing point is the preimage of the fixed point~0 and is
  also a critical point. But when the preimage of~0 is
  solved for numerically, it is sometimes slightly negative and
  sometimes slightly positive. If the final lap had an endpoint at the 
  critical point, this would lead to a failure to find the corresponding
  framing point in step~\ref{normalize}.} 
to the left or right; it is simplest to take them to
be bounded by $\pm\infty$.

As initial values for the marked points, we choose 
    ${x_{1}^{[0]}, x_{2}^{[0]},\ldots,x_{n-1}^{[0]}}$
    to be equally spaced between $x_0=0$ and $x_n=1$,
    and take the initial map $f_0$ to be the piecewise linear map obeying
    the given combinatorics. 

\item \label{mapmaking} \textbf{Mapmaking.}
  Given a map $f_{\ell-1}$, the first step of the iterative process is to
  choose $\breve{f_\ell}$ with the correct critical values: 
  we must determine a map $\breve{f_\ell}$  with critical points 
  $\breve{c_j}^{[\ell]}$
  corresponding to the desired critical values.
  That is,
  \begin{equation}\label{eqMapmake}
    \breve{f_{\ell}}(\breve{c_j}^{[\ell]})=x^{[\ell-1]}_{m_j}=v^{[\ell]}_j, \quad
    \text{where $m_j$ is the index of the image of $\breve{c_j}^{[\ell]}$} \,.
  \end{equation}
  This is done by inverting the map $\Phi$ of \autoref{sec3d-l1}. Given
  the critical-value vector $(v_1, v_2, \ldots, v_r)$, we compute the 
  successive distances 
$$s_1=|v_2-v_1|,\; s_2=|v_3-v_2|,\; \ldots\,,\; s_{r-1}=|v_r-v_{r-1}|,$$
and then use Newton's method to find $\delta_i$ so that 
$$(\delta_1, \ldots, \delta_{r-1}) \approx \Phi^{-1}(s_1,\ldots,s_{r-1}).$$

While Newton's method can be unpredictable,
a good initial choice is to use a scaled version of the critical
points for the corresponding Chebyshev polynomial of the first
kind\footnote{This is convenient since Chebyshev polynomials,
like cubic polynomials, have only two critical values; so that
the distance $s_j$ between consecutive critical values is always constant.}.
Specifically, take   \[\rho_j = \frac{2}{4^{1/r}}\Bigl(
\cos\bigl( j\pi/r \bigr) - \cos\bigl( (j-1)\pi/r \bigr)     \Bigr) \]
as initial point for the Newton iteration
$
  \vec\rho\, \mapsto\, 
  \vec\rho - \bigl(\Phi'(\vec{\rho})\bigr)^{-1}\,\Phi(\vec\rho)
$;
see \autoref{conjNewton}.
This yields a map $\breve{f_\ell}$ with the desired critical values
to within any given tolerance, although the critical points $\breve{c_j}$
of $\breve{f_\ell}$ will not necessarily be in $[0,1]$.

\item \label{normalize} \textbf{Normalization.}
  In order to have a map $f_\ell \in K_\R$, we first need to determine
  the framing points by solving the appropriate framing condition.
  We find $A$ in the initial lap and $B$ in the final lap such that
  $\breve{f_\ell}(A)$ and $\breve{f_\ell}(B)$ are either $0$ or
  $1$ as given by the combinatorics.
  Then we let $f_\ell = \breve{f_\ell}\circ \mu$ and
  $c_j=\mu(\breve{c_j})$, where $\mu$ is the appropriate linear map
  with $\mu(A)=0$, $\mu(B)=1$.

  As the degree grows, numerical uncertainty in the framing points can
  cause the resulting map $f_\ell$ to fail to achieve the desired
  numerical tolerances, in which case we repeat~\ref{mapmaking} with
  increased precision, starting from the~$c_j$ after normalization,
  then normalize the resulting map again.

\begin{figure}[!htb]
\centerline{\includegraphics[height=.4\textwidth]{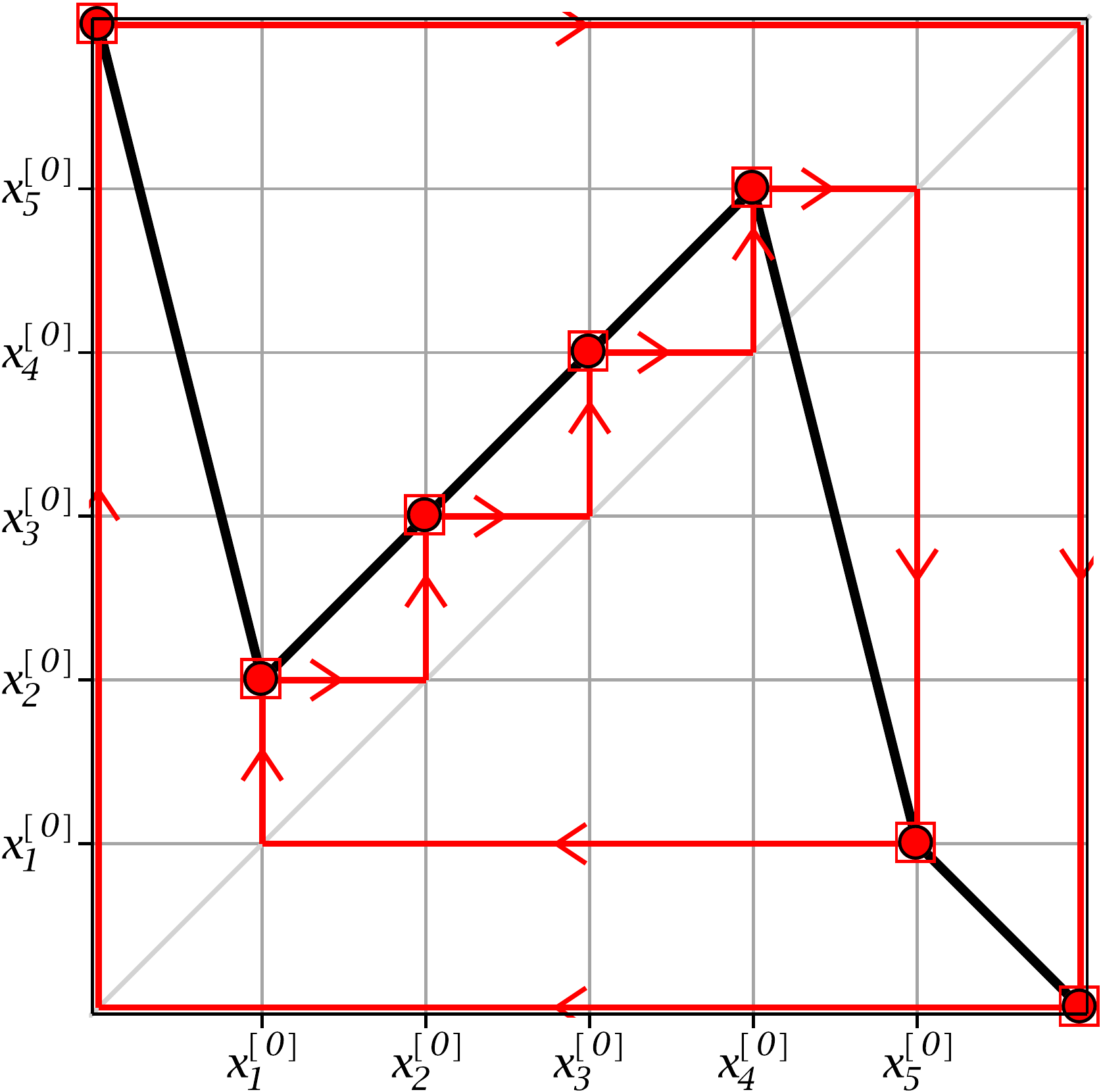} \hfil
  \includegraphics[height=.4\textwidth]{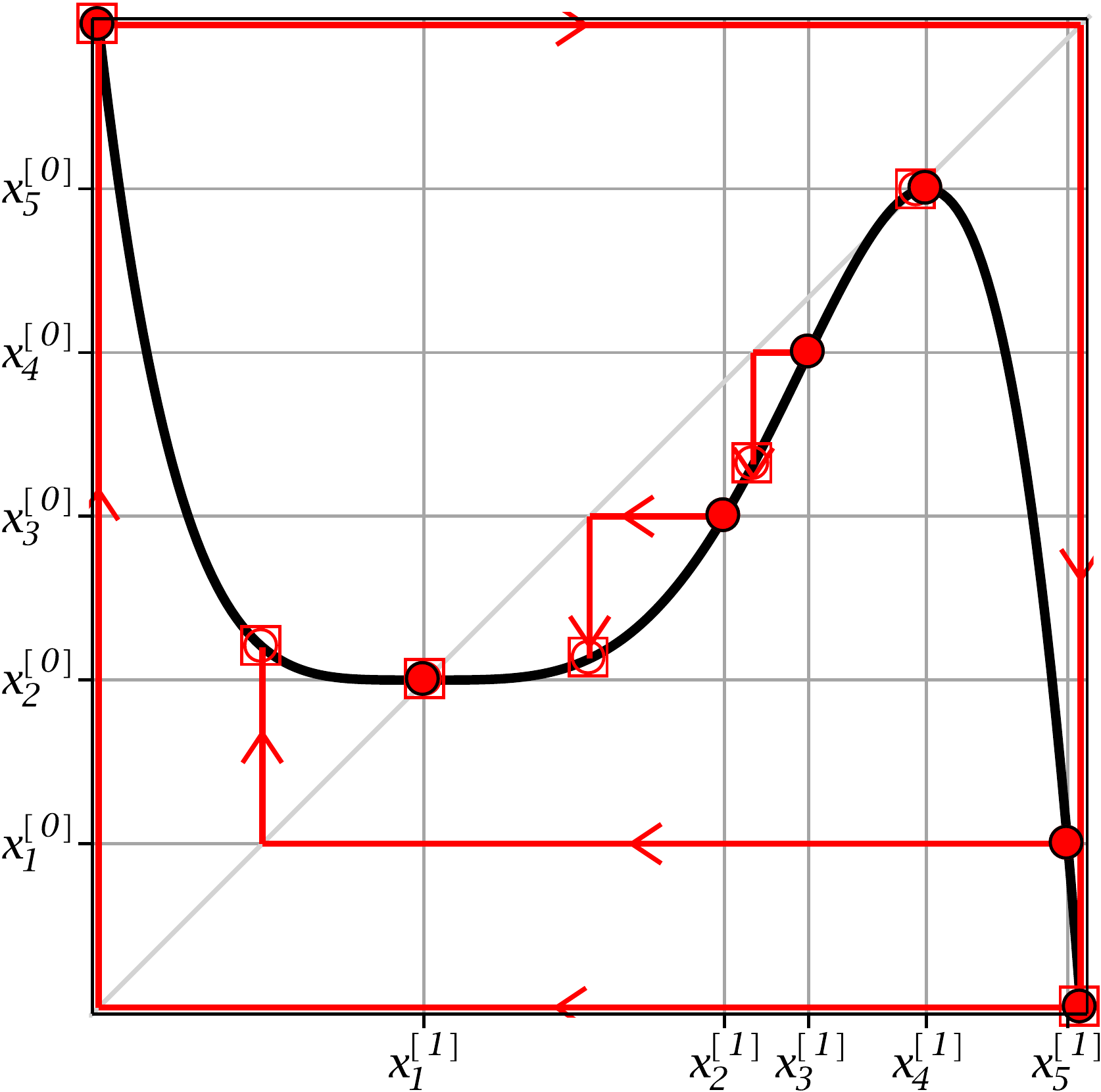}
}
\caption{\label{f-6234510-a} \capf
The piecewise linear map $f_0$ (on the left) and the first step $f_1$
(on right) for combinatorics $(6,2^4,3,4,5,1,0)$.
Here (and in \autoref{f-6234510-b}), the image of a marked point is
indicated by an open square on the graph. 
When passing from the piecewise linear $f_0$ to the polynomial  $f_1$,
observe that while $f_1(x_k^{[1]})=f_0(x_k^{[0]})$, 
the points $x_k$ in the domain have all moved significantly.  Further,
at step~1, the dynamics are significantly off (the image of each
$x_k^{[1]}, 1\le k < 5$ is too far to the left).
}
\end{figure}

\item \label{pullback} \textbf{Pullback.} 
  We now have the unique map $f_\ell \in K_\R$ which satisfies
  Equation~\eqref{eqMapmake};
  this is the map $f_\ell$ on $[0,1]$ with the specified critical values
  $v_j^{[\ell]}$ and known critical points $c_j^{[\ell]}$.
  Since each of these points $c_j^{[\ell]}$ is a point $x_{m_j}^{[\ell]}$,
  we just need to find the remaining noncritical $x_k^{[\ell]}$.
  For each of these, we solve (numerically) 
  $$f_\ell(x_k^{[\ell]})=x_{m_k}^{[\ell-1]},$$
  where $m_k$ is the index of the image of $x_k$ given by the combinatorics and
  $x_k^{[\ell]}$ lying in the appropriate lap, as discussed in \autoref{s3}.
  See Figures~\ref{f-6234510-a} and \ref{f-6234510-b}.

\item \label{howclose} \textbf{Goodness of fit.} 
  To measure how well $f_\ell$ satisfies the conditions, calculate the
  distance~$\epsilon_\ell$ between the images 
 $f_\ell(x_k^{[\ell]})$ and the corresponding points in 
  the orbit $y_k^{[\ell]}  = x_{m_k}^{[\ell]}$, that is,
  \[ \epsilon_\ell =
     \frac{1}{n}
     \Bigl( { \sum_{j=0}^{n}
          \left.  f_\ell(x_k^{[\ell]})^2 - (y_k^{[\ell]})^2\right.}
     \Bigr)^{1/2} \,. \]
  If $\epsilon_\ell$ is not sufficiently small, increment $\ell$ and
  repeat steps~\ref{mapmaking}--\ref{howclose}.
\end{enumerate}

\begin{figure}[thb]
\centerline{\includegraphics[height=.4\textwidth]{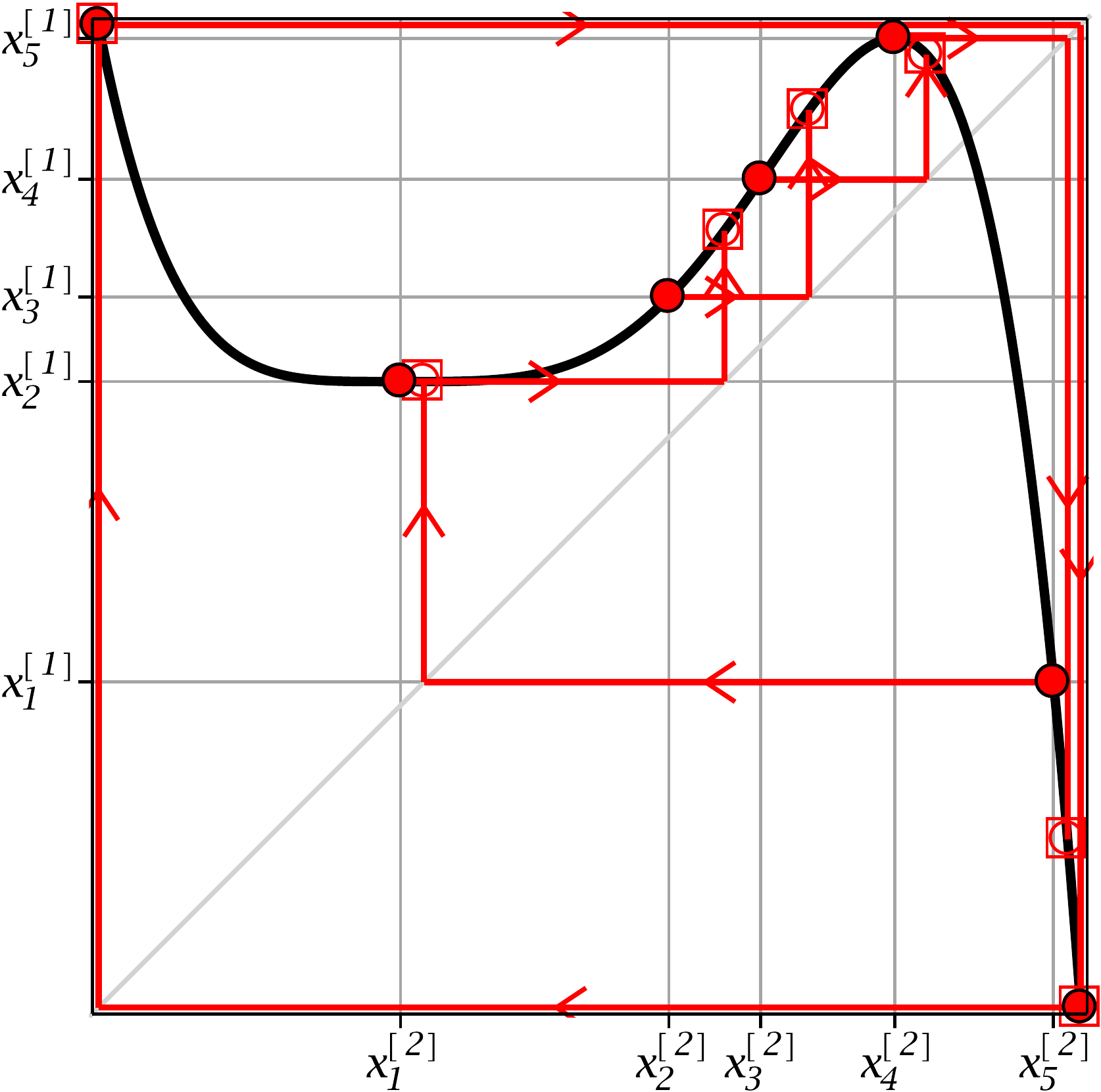} \hfil
  \includegraphics[height=.4\textwidth]{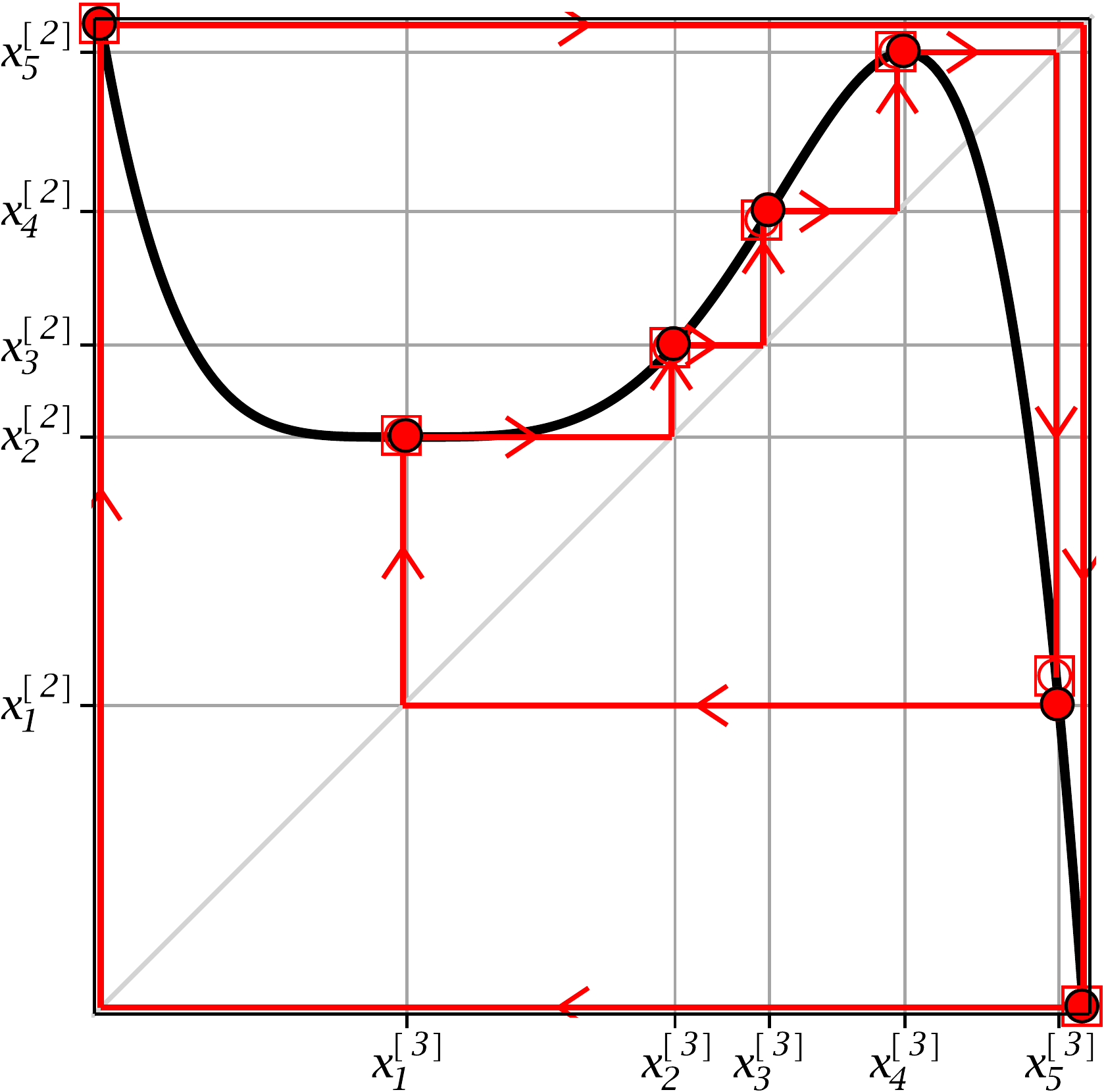}
}
\caption{\label{f-6234510-b} \capf
Steps~2 and 3 of the method for the combinatorics of
\autoref{f-6234510-a}. Note that Step~3 (on the right) is becoming
close to the desired limit dynamics (although only at step~11 is
the error within $10^{-5}$).
}
\end{figure} 


As noted earlier, the global behavior of Newton's method can be
unpredictable, even for diffeomorphisms.  Hence, to determine a
polynomial with the given critical values in step~\ref{mapmaking}, we
need to select an appropriate initial approximation.

Since Chebyshev polynomials have distinct critical points and only two
critical values, we can easily use them to construct polynomials
with critical points $c_j$ for which
$\Phi(\rho_1, \ldots, \rho_{r-1})=(1,1,\ldots,1)$ (with
$|c_j-c_{j+1}|=\rho_j$). 
Assuming the following conjecture, Newton's method will always 
converge to the desired solution when started from such an initial
condition. 

\begin{conj}\label{conjNewton}
Given $(s_1,\ldots, s_{r-1})$ with $0 < s_j\leq 1$ for all $j$, 
suppose that $(\rho_1,\ldots, \rho_{r-1})$ is chosen so that
$\Phi(\rho_1, \ldots, \rho_{r-1}) = (\sigma_1, \ldots,
\sigma_{r-1})$, with $s_j \le \sigma_j\le 1$ for all $j$. Then  Newton's method,
when started with $(\rho_1, \ldots, \rho_{r-1})$ as the initial condition,
will converge to the desired solution $(\delta_1,\ldots,\delta_{r-1})
=\Phi^{-1}(s_1,\ldots, s_{r-1})$. 
\end{conj}

While we haven't quite been able to establish \autoref{conjNewton}, 
observe that all entries of the first derivative matrix
$\Phi'$ are non-negative in this region,
as are the partial derivatives of each entry. 
These properties certainly simplify the situation, and in most cases
Newton's method decreases monotonically in each coordinate towards
$(\delta_1,\ldots,\delta_{r-1})$.

Even without assuming the conjecture, one can apply a modified version
of Newton's method to ensure convergence (see \cite[ch.3]{Deuflhard}, 
for example).  We have not found a situation in which this was necessary. 

\bsk\bsk
\section{Further Examples: The Non-Expansive Case}\label{s4}

 In this section, after presenting another typical example, we move to
examples with non-expansive combinatorics. In the non-expansive case,
the Thurston algorithm converges 
only in a much weaker sense, as follows from a theorem of Selinger~\cite{S}.
\medskip

\begin{figure}[!htb]
\centerline{\includegraphics[height=.35\textwidth]{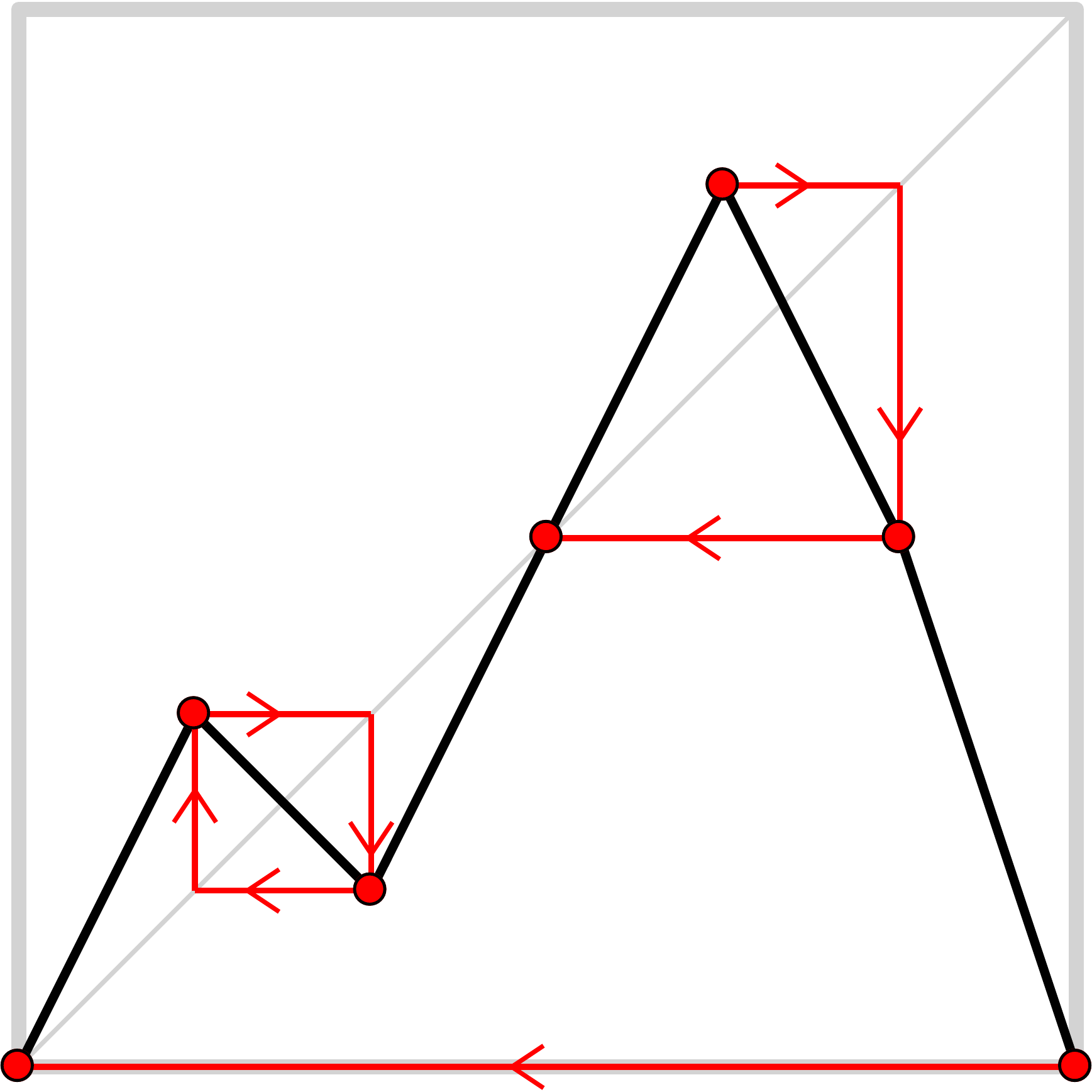} \hfil
  \includegraphics[height=.35\textwidth]{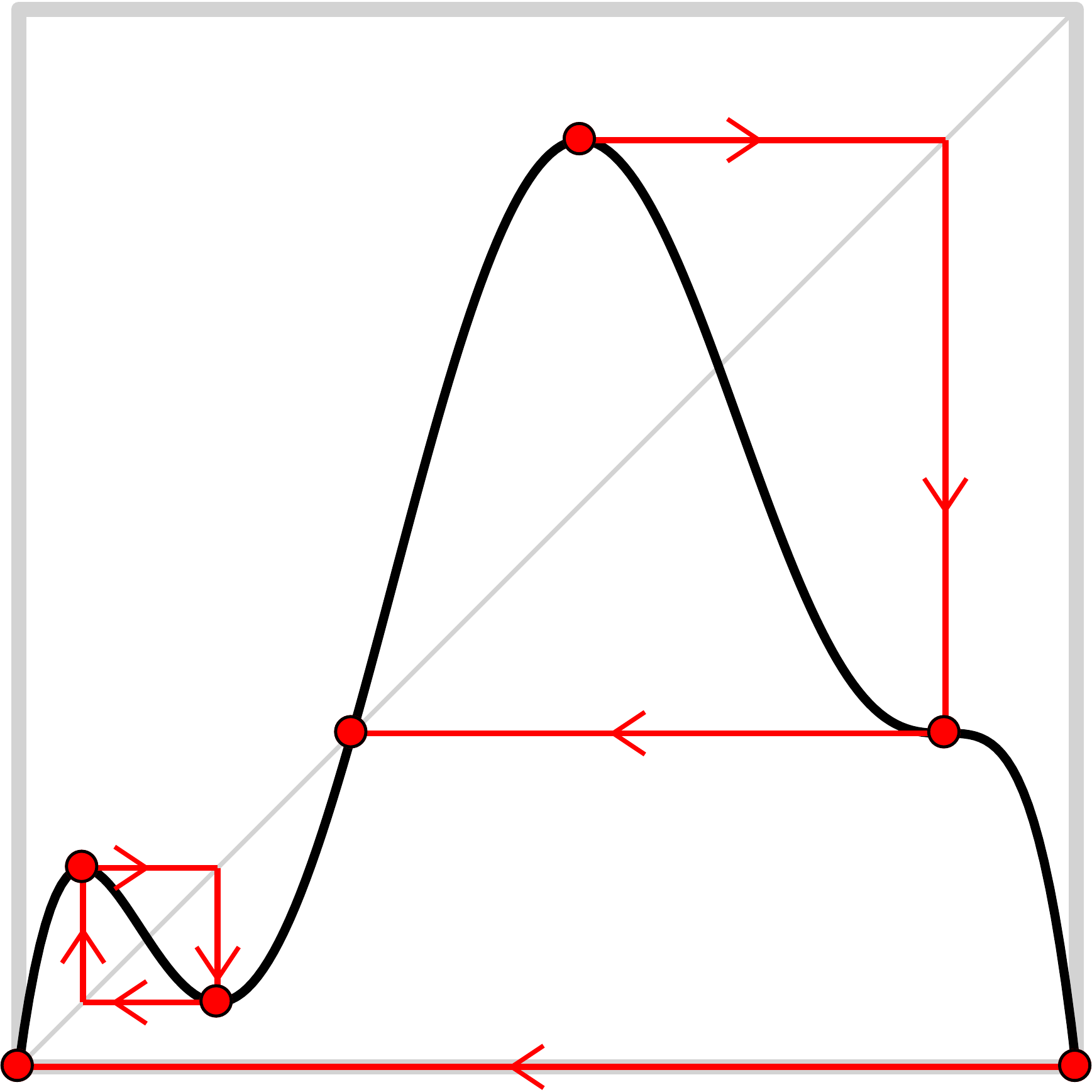}%
}
\caption{\label{f-ex1} 
{\capf On the left is the PL model for the combinatorics
  $(0,2,1,3,5,3^3,0)$, with
  \notAMS{mapping pattern }
  $~~\du{x_1}\leftrightarrow\du{x_2}$ and
  $\du{x_4}\mapsto \underline{\du{x_5}}\mapsto x_3\mapstoself$%
\AMSonly{ as a mapping pattern}.
 On the right is the corresponding polynomial map of degree six.}}\bsk
\end{figure}

\begin{figure}[!htb]
\centerline{\includegraphics[height=.35\textwidth,align=t]{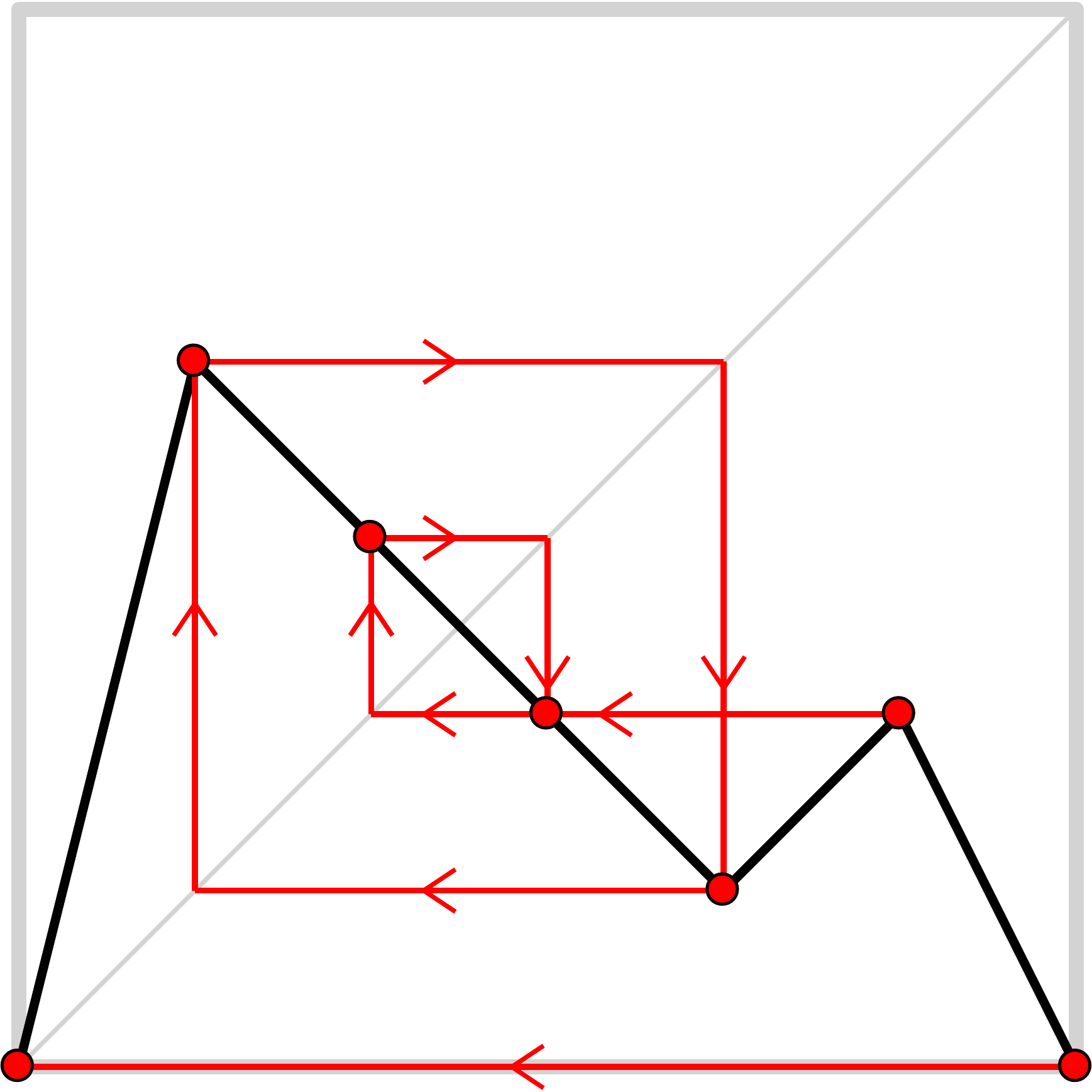}\hfil
  \includegraphics[height=.35\textwidth,align=t]{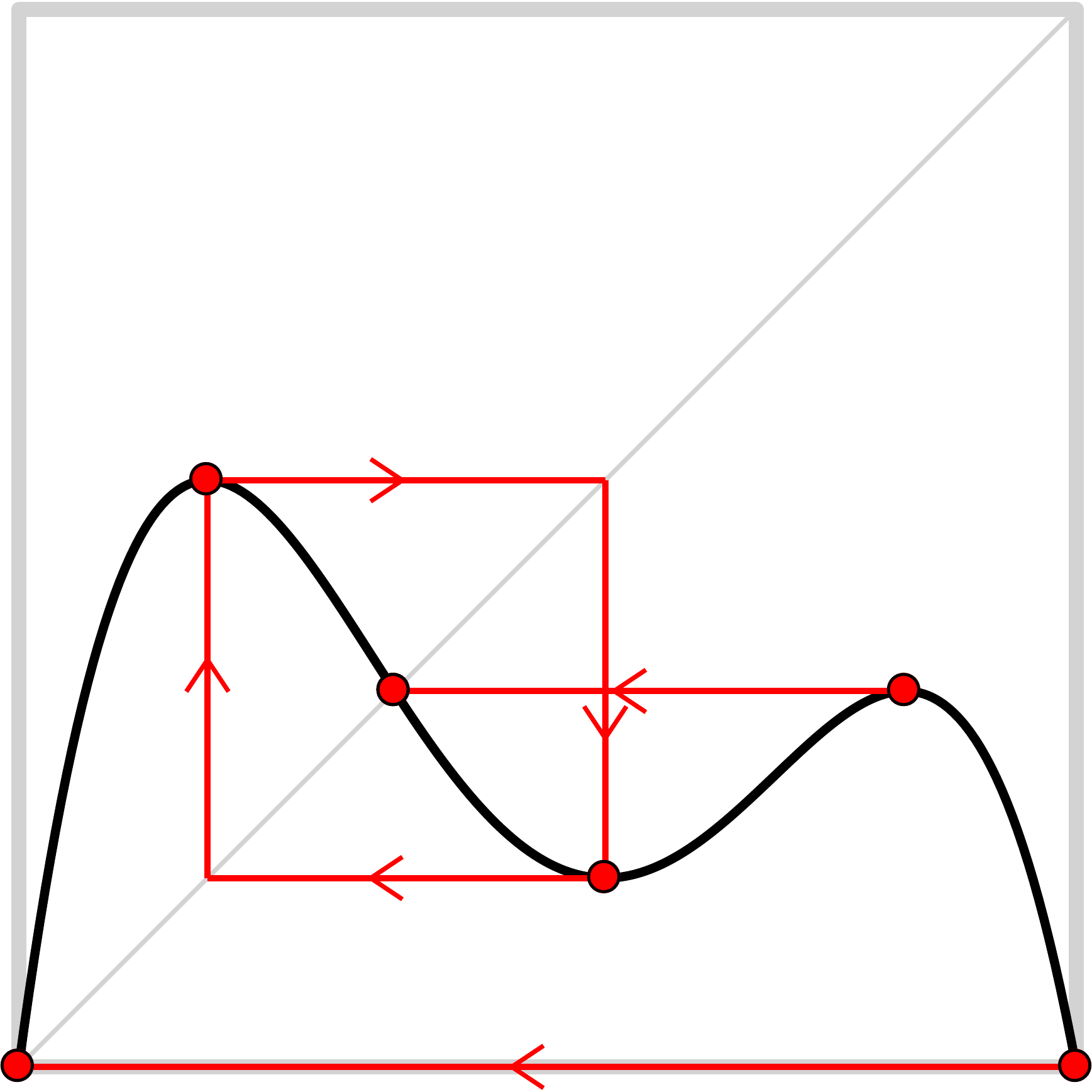}
}

\caption{\label{f-ex4} {\capf On the left is the
PL model for combinatorics $(0,4,3,2,1,2,0)$ with mapping pattern
 $\du x_1 \leftrightarrow \du{x_4}~,\quad \du{x_5}\mapsto x_3\leftrightarrow x_2$. On the right
is the corresponding polynomial map of degree four.
 Here the interval $[2,3]$ of the PL model on the left
   has shrunk to the circled fixed point, and the combinatorics
   has simplified to $(0,3, 2,1,2,0)$, with mapping pattern 
 $~~\du{x_1}\leftrightarrow \du{x_3}~,\quad \du{x_4}\mapsto x_2 \mapstoself$.}}
\end{figure}

Figure \ref{f-ex1} shows an example which satisfies all six of the
conditions
of \autoref{s1}.  For such examples, all of the marked points will  remain distinct 
in the limit. On the other hand, \autoref{f-ex4} shows an example of 
 combinatorics which violates the expansiveness condition \ref{ComboExpansive}~of
 \autoref{s1}: no forward image of the interval 
$[2,3]$ contains a critical point. Applying the Thurston algorithm
to this example, the interval $[2,3]$ shrinks to a point, as shown in the graph
to the right in \autoref{f-ex4}, so that the combinatorics becomes simpler.
\bsk

\begin{figure}[!htb]
\centerline{\includegraphics[width=.35\textwidth]{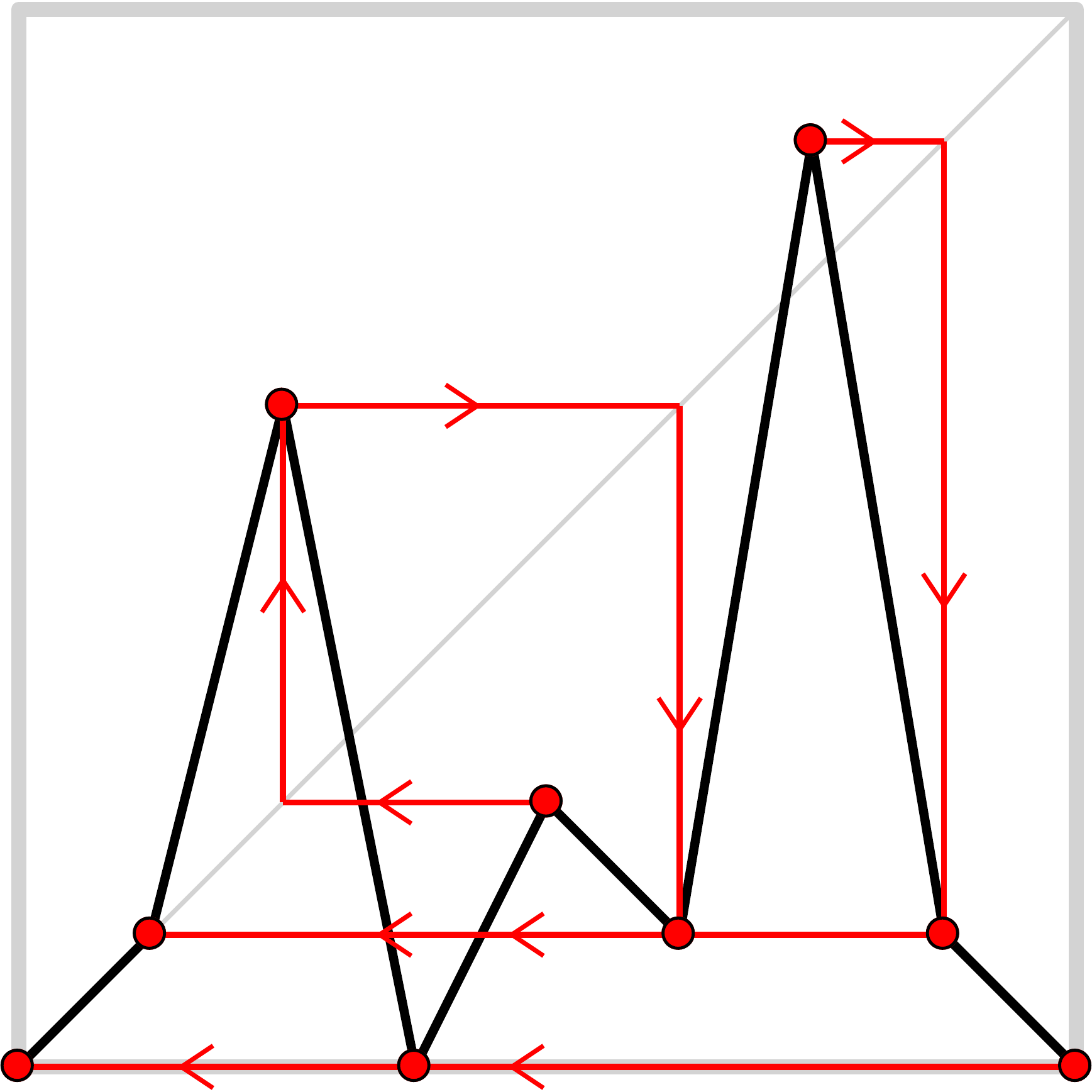} \hfil
  \includegraphics[width=.35\textwidth]{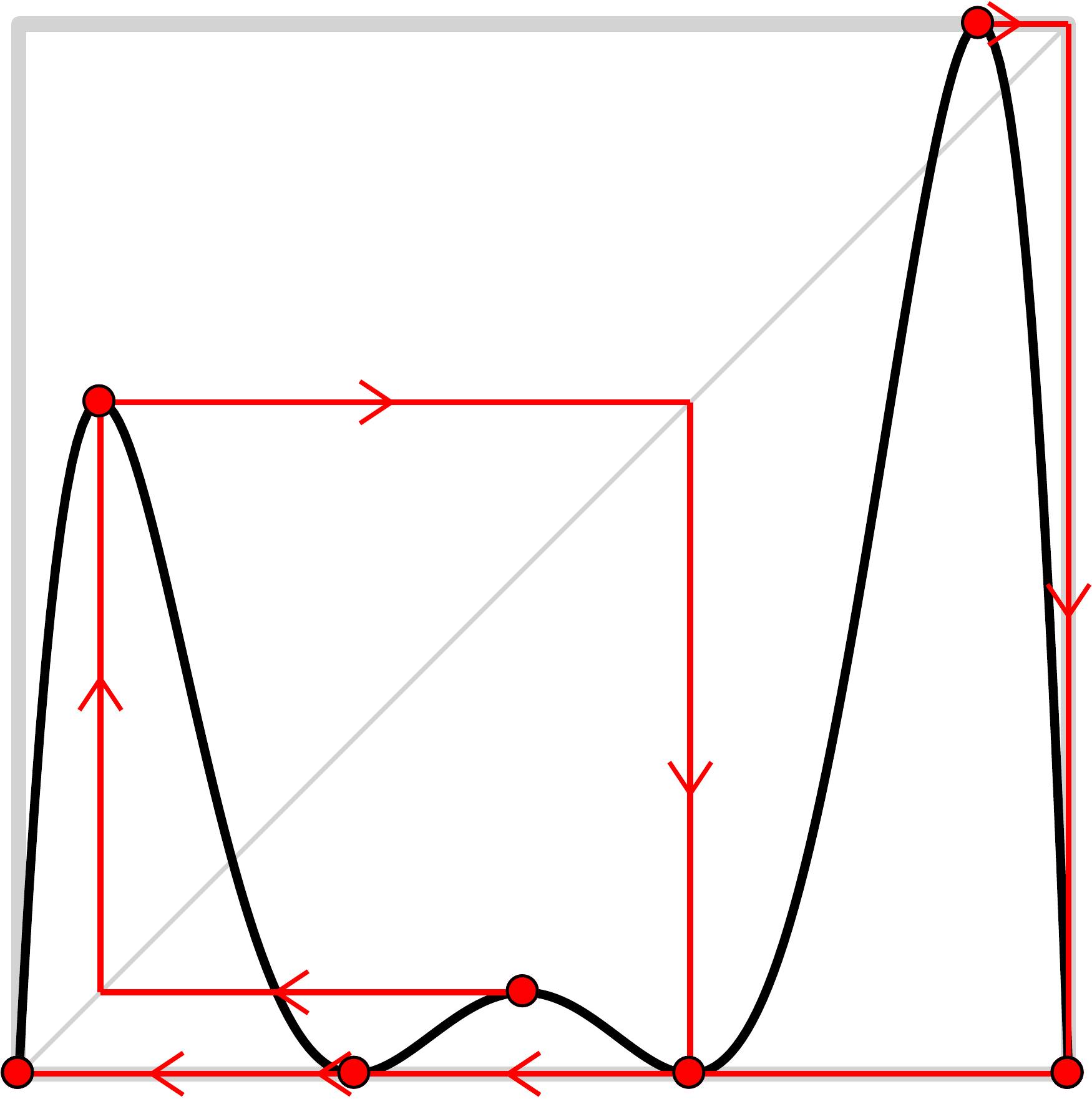}
}
  \caption{\label{f-ex5} \capf On the left: a piecewise linear map with
    combinatorics  $(0,1,5,0,2,1,7,1,0)$. Here every critical orbit ends 
    either at the fixed point zero or at the fixed point with number one.
    However  the initial interval $[0,1]$ and the
     final interval $[7,8]$ fail to be expansive; so on the right the
     associated polynomial map is simpler,
     with combinatorics $(0,4,0,1,0,6,0)$. 
    The mapping pattern in the limit is
     {$\du{x_3}\mapsto\du{x_1}\mapsto\du{x_4}\mapsto x_0 \mapstoself$,}~~
     {$\du{x_5}\mapsto x_6\mapsto x_0 \mapstoself$,}
      and $~~\du{x_2}\mapsto x_0\mapstoself$,~~ with all critical orbits 
 ending at the repelling fixed point $x_0$.}
 \end{figure}

\autoref{f-ex5} shows a similar example. Here both of the intervals $[0,1]$ and
$[7,8]$ collapse to points under the Thurston algorithm, so that again the final 
   polynomial map is simpler than the original PL map.

 These examples illustrate the following statement, which follows relatively
  easily from the work of Nikita Selinger \cite[Proposition 6.2]{S}.
  We want to thank Kevin Pilgrim, Thomas Sharland,  as well as Selinger himself,
  for  pointing this out to us.
\ssk

\begin{asse}
  Even if the given combinatorics does not satisfy the expansiveness
 condition, the successive polynomial approximations given by the Thurston
 algorithm will still converge locally uniformly to a critically finite
 polynomial. However this limit polynomial will have simpler combinatorics.
 More precisely, every edge of
 the piecewise linear model which is not expansive will collapse to a point.
 \end{asse}

 We will not attempt to provide further details, but encourage the interested
 reader to study Selinger's paper. 


 \section{Coefficients.}\label{s5}
\AMSorNot{}{\vspace*{-.5\baselineskip}} 
%
\begingroup \renewcommand{\arraystretch}{1.5}
  \newcommand*{\thead}[1]{\multicolumn{1}{c|}{\textsf{#1}}}
  \begin{center} 
  \small
  \begin{longtable}{|c|p{.68\textwidth} |c|c|}
  \caption{\label{coeff-tab} \capf
    The coefficients for the polynomials in several of the figures, and 
    the corresponding error estimates, evaluated as 
    $\tfrac{1}{n}\sqrt{ {\sum (f(x_i)-y_i)^2}}$ where the $x_i$ are the 
    marked points, and $y_i$ is the desired value for $f(x_i)$.    
    The error estimates marked with * represent
    initial steps for the Thurston algorithm, while the remaining
    estimates represent estimates after many steps of the algorithm. Note
    that the coefficients for \autoref{f-Thalg} are beginning to
    converge towards the coefficients for \autoref{f-thex-lim}. Similarly, the
    coefficients for Figures~\ref{f-6234510-a}(r) and \ref{f-6234510-b}
  are beginning to converge towards the Thurston limit.}
\\
\hline
\textsf{Fig.}&\thead{polynomial} & \thead{error bound} &\thead{Iter} \\
\hline\hline
\endfirsthead
 \captionsetup{singlelinecheck=off}
   \caption[]{\textit{Continued from previous page}}\\
\hline
\textsf{Fig.}&\thead{polynomial} & \thead{error bound} &\thead{Iter} \\
\hline\hline
\endhead
  \hline
  \multicolumn{4}{r}{\textit{\sffamily \autoref{coeff-tab} continues on next page}}\\
\endfoot 
  \hline
\endlastfoot

\ref{f2} 
&$7.121692805x-17.64597623x^2+11.52428342x^3$
&$1.24 \times 10^{-8}$ & 20\\
\ref{f-Thalg}(l)
&$15.332055x -92.795911x^2 +225.00679x^3 -242.71367x^4 +96.170733x^5$
&$0.037^*$  & 1 \\  
\ref{f-Thalg}(r)
&$18.069912x -112.83091x^2 +273.38011x^3 -292.41971x^4 +114.80059x^5$
& $0.0038^*$  & 2 \\
\ref{f-thex-lim} 
&$18.163069x -113.72167x^2 +276.22221x^3 -296.09149x^4 + 116.42789x^5$
&$1.84\times 10^{-6}$ & 13 \\ 
\ref{f-03214}(l)
&$6x -15x^2+10x^3$
&$1\times 10^{-13}$ & 14\\ 
\ref{f-03214}(r)
&$0.20557075x-181.7478872x^2+855.1404749x^3-2244.547436x^4+3255.216137x^5-2427.230116x^6+723.9632564x^7$
&$4.13\times10^{-8}$ & 18 \\
\ref{f-6234510-a}(r)
&$1 -8.73730x +44.7494x^2 -110.928x^3 +130.960x^4 -57.0449x^5$
&$0.0791^*$ & 1 \\
\ref{f-6234510-b}(l)
&$1-5.10905x+28.0816x^2-74.57010x^3+93.9995x^4-43.4011x^5$ 
&$0.0144^*$ & 2 \\
\ref{f-6234510-b}(r)
&$1.-5.82395x+31.5803x^2-82.7518x^3 +102.974x^4 -46.9785x^5$
&$0.0021^*$ & 3\\
\ref{f-ex1}(r) 
&$7.494214522x-97.01797994x^2+457.9211574x^3-913.0123135x^4+811.6279094x^5-267.0129879x^6$   
& $3.54 \times 10^{-9}$ &25\\
\ref{f-ex4}(r) 
&$7.45977893x-32.0733758x^2+47.0904007x^3-22.4768041x^4$
&$5.49 \times 10^{-9}$ &36 \\
%
\ref{f-ex5}(r)
&$20.15184092x-208.9317665x^2+827.5262978x^3-1559.747539x^4+1400.650082x^5-479.6489149x^6$
&$6.34 \times 10^{-8}$ & 12\\                                                        %
\hline
\end{longtable}
  \end{center}
  \endgroup


\end{document}